\documentclass[12pt]{article}

\usepackage{amssymb}
\usepackage{amsmath}
\usepackage{color}
\usepackage{graphics}
\usepackage{epsfig}
\usepackage{hyperref}

\newtheorem{claim}{\bf \t}[part]

\newtheorem{Lemma}{Lemma}[part]

\newtheorem{Remark}{Remark}[part]
\newtheorem{Theorem}{Theorem}[part]

\numberwithin{Assumption}{section} \numberwithin{Corollary}{section}
\numberwithin{Definition}{section} \numberwithin{equation}{section}
\numberwithin{Example}{section} \numberwithin{Lemma}{section}
\numberwithin{Proposition}{section} \numberwithin{Remark}{section}
\numberwithin{Theorem}{section}

\def\v{\varepsilon}

\def\x{\xi}
\def\t{\theta}
\def\T{\Theta}

\def\a{\alpha}
\def\b{\beta}
\def\g{\gamma}
\def\d{\delta}
\def\l{\lambda}
\def\f{\frac}

\def\r{\rho}

\def\s{\sigma}

\def\o{\omega}
\def\di{\displaystyle}
\def\i{\infty}

\def\text#1{{\rm #1}}
\begin{document}
\date{}
\title{\Large \bf Global classical solutions to the two-dimensional  compressible Navier-Stokes equations in $\mathbb{R}^2$ }
\author{\small \textbf{Quansen Jiu},$^{1,3}$\thanks{The research is partially
supported by National Natural Sciences Foundation of China (No.
11171229) and Project of Beijing Education Committee. E-mail:
jiuqs@mail.cnu.edu.cn}\quad
  \textbf{Yi Wang}$^{2,3}$\thanks{The research is partially supported by National Natural Sciences Foundation of China (No.
11171326) and by the National Center for Mathematics and
Interdisciplinary Sciences, CAS. E-mail: wangyi@amss.ac.cn.}\quad
and \textbf{Zhouping Xin}$^{3}$\thanks{The research is partially
supported by Zheng Ge Ru Funds, Hong Kong RGC Earmarked Research
Grant CUHK4042/08P and CUHK4041/11P, and a grant from the Croucher
Foundation. Email: zpxin@ims.cuhk.edu.hk}} \maketitle \small $^1$
School of Mathematical Sciences, Capital Normal University, Beijing
100048, P. R. China

\small $^2$ Institute of Applied Mathematics, AMSS, CAS, Beijing 100190, P. R. China

\small $^3$The Institute of Mathematical Sciences, Chinese
University of HongKong, HongKong\\

 {\bf Abstract:} In this paper, we prove the global well-posedness of the classical solution to the  2D Cauchy problem of
 the compressible Navier-Stokes equations   with arbitrarily large initial data when the shear viscosity $\mu$ is
 a positive constant and the bulk viscosity $\l(\r)=\r^\b$ with $\b>\frac43$. Here the initial density keeps a
 non-vacuum states $\bar\rho>0$ at far fields and our results generalize the ones by Vaigant-Kazhikhov \cite{Kazhikhov}
 for the periodic problem and by Jiu-Wang-Xin \cite{JWX3} and Huang-Li \cite {hl2} for the Cauchy problem with
 vacuum states $\bar\rho=0$ at far fields. It shows that the solution will not develop the vacuum states in any finite time
 provided the initial density is uniformly away from vacuum. And the results also hold true when the initial data
  contains vacuum states in a subset of $\mathbb{R}^2$ and the natural
 compatibility conditions are satisfied. Some new weighted estimates are obtained  to establish the
 upper bound of the density.

{\bf Key Words:} compressible Navier-Stokes equations, Cauchy
problem, global well-posedness, large data, vacuum


\section{Introduction } \setcounter{equation}{0}
\setcounter{Assumption}{0} \setcounter{Theorem}{0}
\setcounter{Proposition}{0} \setcounter{Corollary}{0}
\setcounter{Lemma}{0} We consider the following
compressible and isentropic Navier-Stokes equations
\begin{eqnarray}\label{CNS}
\left\{ \begin{array}{ll}
\partial_t\rho+{\rm div}(\rho u)=0,  \\
\partial_t(\rho u) + {\rm div} (\rho u\otimes u) + \nabla P(\r)=\mu\Delta u+\nabla((\mu+\l(\r)){\rm div}u), &\quad x\in\mathbb{R}^2, t>0,
\end{array}
\right.
\end{eqnarray}
where $\rho (t, x)\geq 0 $, $u(t, x)=(u_1,u_2)(t,x) $ represent the
density and the velocity of the fluid, respectively. And
$x=(x_1,x_2)\in\mathbb{R}^2$ and
$t\in[0,T]$ for any fixed $T>0$. Here it is assumed that the shear viscosity $\mu>0$ is a positive constant and the bulk viscosity
\begin{eqnarray}\label{1.2}
\l(\r)=\r^\b
\end{eqnarray}
with $\b>0$  in general such that the operator
$$
\mathcal{L}_\r u\equiv\mu\Delta u+\nabla((\mu+\l(\r)){\rm div}u)
$$
is strictly elliptic. The pressure function is given by
$P(\r)=A \r^{\gamma },$
where $\gamma >1$ denotes the adiabatic exponent and $A>0$ is a
constant which is normalized to be $1$ for simplicity. We impose the initial values as
\begin{equation}\label{initial-v}
(\r,u)(t=0,x)=(\r_0,u_0)(x)\rightarrow(\bar\r,0),\quad {\rm as}\quad
|x|\rightarrow+\i,
\end{equation}
where $\bar\r>0$ is a given positive constant.

In the case that both the shear and bulk
viscosities are positive constants, there are a large number of literatures  on the well-posedness theories of the
compressible Navier-Stokes equations.  In particular, the
one-dimensional theory is rather satisfactory, see
\cite{Hoff-Smoller, lius, Ka, KS} and the references therein. In
multi-dimensional case,
 the local well-posedness theory of classical solutions  was
 established   in the absence of
 vacuum (see \cite{Nash},  \cite{Itaya} and  \cite{Tani}) and the global well-posedness theory of classical solutions  was
 obtained  for initial data
  close to a non-vacuum steady state (see \cite{MN},  \cite{Hoff1},   \cite{Dachin},  \cite{CMZ} and references therein).  The local well-posedness
 of   classical solutions containing vacuum was studied  by Cho-Kim \cite{CK1}
 and Luo \cite{Luo} and the global well-posedness
 of classical solutions to the 3D isentropic compressible Navier-Stokes equations with small energy  was proved by Huang-Li-Xin \cite{hlx4}.
    For  the large  initial data permitting vacuums,  the global existence of weak solutions was
    investigated in  \cite{Lions}, \cite{F1}, \cite{JZ}.
   It should be noted that if the initial data are arbitrarily large and the vacuums are permitted,  the solution will also contain possible vacuums and one could not expect the global well-posedness in general, see  \cite{Xin}  \cite{R} and
   \cite{Xin-Yan} for blow-up results of  classical solutions.

The case that both the shear and bulk viscosities  depend on the
density has also received a lot attention recently, see \cite{BD1,
BDL, Dachin, DWZ, GJX, JXZ, JZ, JWX, JWX1, JX, LLX, LXY, MV, YYZ,
YZ2, ZF} and the references therein. When deriving by Chapman-Enskog
expansions from the Boltzmann equation, the viscosity of the
compressible Navier-Stokes equations  depends on the temperature and
thus on the density for isentropic flows. Moreover, in  geophysical
flows, the viscous Saint-Venant system for the shallow water
corresponds exactly to a kind of compressible Navier-Stokes
equations with density-dependent viscosities. However, except for
the one-dimensional problems, few results are available for the
multi-dimensional problems and  even the short time well-posedness
of  classical solutions in the presence vacuum remains open.

The system  \eqref{CNS} was first proposed and studied
by Vaigant-Kazhikhov in \cite{Kazhikhov}. For the periodic problem on the torus $\mathbb{T}^2$ and under assumptions that  the initial density is
uniformly away from vacuum and $\beta>3$ in \eqref{1.2},  Vaigant-Kazhikhov   established the well-posedness of the classical solution to \eqref{CNS} in \cite{Kazhikhov} and the global existence and large time behavior of weak solutions was stuided by   Perepelitsa in \cite{P}.
Recently, Jiu-Wang-Xin \cite{JWX2} improved the result in
\cite{Kazhikhov} and obtained the global well-posedness of the
classical solution to the periodic problem with large initial data
permitting vacuum. Later on, Huang-Li relaxed the index $\beta$ to be $\beta>\frac
43$ and studied the large time behavior of the solutions in
\cite{hl1}. However, all the above results are concerned with the 2D
periodic problems. For the 2D Cauchy problems with  vacuum states at far fields, Jiu-Wang-Xin \cite{JWX3} and Huang-Li \cite{hl2} independently considered the global well-posedness of classical solution in different weighted spaces.

 In the present paper, we study the global
well-posedness of the classical solution to the Cauchy problem
\eqref{CNS}-\eqref{initial-v} with large data which keeps a
non-vacuum states $\bar\rho>0$ at  far fields. In particular, our
results show that the solution will not develop the vacuum states in
any finite time provided the initial density is uniformly away from
vacuum. The results of this paper  generalize the ones by
Vaigant-Kazhikhov in \cite{Kazhikhov} to  the Cauchy problem and the
index $\beta$ is relaxed to be $\beta>\frac43$. The results also
improve ones by Jiu-Wang-Xin \cite{JWX3} and Huang-Li \cite {hl2}
for the Cauchy problem with vacuum states $\bar\rho=0$ at far
fields. Moreover, the results  hold true if the initial data
contains vacuum states in a subset of $\mathbb{R}^2$ under
appropriate compatibility conditions (see \eqref{cc} in Theorem
1.2).

 To study the global well-posedness of the classical solution of the compressible Navier-Stokes equations,
 it is crucial to obtain the uniformly upper bound of the density. To do that, similar to  \cite{Kazhikhov}, \cite{JWX2} and \cite{JWX3}, we first obtain any $L^p (2\le p<\infty)$ estimates of the density $\rho-\bar\rho$  and then obtain the estimates of the first order derivative of the velocity. A new transport equation \eqref{transport-e} is derived by means of  the effective viscous flux $F=(2\mu+\l(\r)){\rm div}u-(P(\r)-P(\bar\r))$ and two new functions $\xi$ and $\eta$ satisfying the elliptic problems
\begin{equation}\label{xi-10}
-\Delta \x={\rm div}(\r u),\qquad
-\Delta \eta={\rm div}[{\rm div}(\r u\otimes u)],
\end{equation}
respectively, which was introduced in \cite{Kazhikhov}. Comparing
with the periodic problem and  the Cauchy problem with vacuum at far
fields, new difficulties will be encountered. Since no integrability
is expected for the density itself, we will decompose the elliptic
problem \eqref{xi-10} into the following two parts:
 \begin{equation}\label{xi-11}
-\Delta \x_1={\rm div}(\sqrt\r u(\sqrt\r-\sqrt{\bar\r})),
\end{equation}
\begin{equation}\label{xi2-12}
-\Delta \x_2=\sqrt{\bar\r} ~{\rm div}(\sqrt\r u).
\end{equation}
 For the elliptic problem \eqref{xi-11}, one can make use  of the similar properties as the periodic case  and the
 Cauchy problem with vanishing density at the far fields thanks to the expected integrability of $\sqrt\r-\sqrt{\bar\r}$.
 For the second elliptic problem \eqref{xi2-12}, since it is expected that $\rho\in L^\infty$ and
  $\sqrt\r u\in L^\i([0,T];L^2(\mathbb{R}^2))$  by the elementary energy estimate, it follows from \eqref{xi2-12}
  that $\nabla\xi_2\in D^1(\mathbb{R}^2)$ which is a homogeneous and critical Sobolev space.
   Therefore, the integrability  of $\xi_2$ can not be derived in a direct way. However, the integrability
    of $\xi_2$ is crucial  to obtain the $L^p (2\le p<\infty)$ estimates of $\r-\bar\r$ and the upper bound of
    the density $\r$. In order to circumvent this difficulty, some new weighted estimates are needed and the
    integrability of the velocity and $\xi_2$ is proved by using Cafferelli-Kohn-Nirenberg type inequality
    \cite{CKN, CW}. It should be remarked  that these  weighted estimates are motivated by our previous
    work \cite{JWX3} and in comparison with  the uniform constant  in \cite{JWX3}, the weight power $\a$  here  depends
    on the ratio $\frac{\l(\bar\r)}{\mu}$.  At the same time, if $\bar\r=0$, then the weight $\a$  is exactly
    same as the one in our previous work  \cite{JWX3}. Moreover, when deriving the first-order derivative estimates
    of the velocity, since $L^p$-integrability $(2\le p<\infty)$ is not available, it would be required to use the $L^\infty$-norm of the density  $\r$
    in a priori way which is motivated by the work \cite{P}. In this way, a $\log$-type inequality of
    the first-order derivative of the velocity can be obtained (see Lemma \ref{lemma-u-der}). Finally, with help of a higher energy estimate in Lemma \ref{lemma-u-d}, one can get a upper bound of the density under the restriction $\beta>\frac 43$ (see \cite {hl1,hl2}).

 Denote the potential energy by
$$
\Psi(\r,\bar\r)=\f{1}{\g-1}\big[\r^\g-\bar\r^\g-\g\bar\r^{\g-1}(\r-\bar\r)\big].
$$

 Our main results can be stated as follows.

\begin{Theorem}\label{theorem1}
Let $\beta>\frac 43$ and $1<\gamma\le 2\beta$. Suppose that the
initial values $(\r_0,u_0)(x)$ satisfy
\begin{equation*}
\begin{array}{ll}
\di 0<c\leq\r_0\leq C, \quad (\r_0-\bar\r, P(\r_0)-P(\bar\r))\in
W^{2,q}(\mathbb{R}^2)\times W^{2,q}(\mathbb{R}^2),\quad u_0(x)\in H^2(\mathbb{R}^2),\\
\di \Psi(\r_0,\bar\r)(1+|x|^{\a})\in L^1(\mathbb{R}^2), \quad
\sqrt{\r_0}u_0 |x|^{\f\a2}\in L^2(\mathbb{R}^2),
\end{array}
\end{equation*}
where $q, c, C$ and $\a$ are positive constants satisfying $q>2$,
$0<c<C$ and
$0<\a^2<\f{4(\sqrt{2+\f{\l(\bar\r)}{\mu}}-1)}{1+\f{\l(\bar\r)}{\mu}}$
respectively.  Then, for any $T>0$, there exists  a unique global
classical solution $(\r,u)(t,x)$ to the Cauchy problem
\eqref{CNS}-\eqref{initial-v} satisfying
$$0<c_1\leq \r\leq C_1$$
for some positive constants $c_1$ and $C_1$. Moreover, one has
\begin{equation}\label{Jan-16-1}
\begin{array}{ll}
\di  (\r-\bar\r,P(\r)-P(\bar\r))(t,x)\in C([0,T]; W^{2,q}(\mathbb{R}^2)),\\
\di \Psi(\r,\bar\r)(1+|x|^{\a})\in C([0,T]; L^1(\mathbb{R}^2)),\quad
\di \sqrt\r u  |x|^{\f\a2}\in C([0,T]; L^2(\mathbb{R}^2)),\\
u\in C([0,T];H^2(\mathbb{R}^2))\cap L^2(0,T;H^3(\mathbb{R}^2)),~~
\sqrt t u\in L^\i(0,T; H^{3}(\mathbb{R}^2)),\\t u\in L^\i(0,T;
W^{3,q}(\mathbb{R}^2)),~~
u_t\in L^2(0,T;H^1(\mathbb{R}^2)),\\
\di \sqrt tu_t\in L^2(0,T; H^2(\mathbb{R}^2))\cap
L^\i(0,T;H^1(\mathbb{R}^2)),~~ t u_t\in L^\i(0,T;
H^2(\mathbb{R}^2)),\\ \sqrt t \sqrt\r u_{tt}\in
L^2(0,T;L^2(\mathbb{R}^2)),~~
 t \sqrt\r u_{tt}\in
L^\i(0,T;L^2(\mathbb{R}^2)),~~t\nabla u_{tt}\in
L^2(0,T;L^2(\mathbb{R}^2)).
\end{array}
\end{equation}
\end{Theorem}

If the initial values contain vacuum states in a subset of
$\mathbb{R}^2$, then the following results can be obtained.

\begin{Theorem}\label{theorem2}
Suppose that the initial values $(\r_0,u_0)(x)$ satisfy
\begin{equation}\label{in-d}
\begin{array}{ll}
\di \r_0\geq0, \quad (\r_0-\bar\r, P(\r_0)-P(\bar\r))\in
W^{2,q}(\mathbb{R}^2)\times W^{2,q}(\mathbb{R}^2),\quad u_0(x)\in H^2(\mathbb{R}^2),\\
\di \Psi(\r_0,\bar\r)(1+|x|^{\a})\in L^1(\mathbb{R}^2), \quad
\sqrt{\r_0}u_0|x|^{\f\a2}\in L^2(\mathbb{R}^2),
\end{array}
\end{equation}
with $q, \a, \g$ and $\b$ being the same as in Theorem
\ref{theorem1}. Suppose that  the  compatibility conditions
\begin{equation}\label{cc}
\mathcal{L}_{\r_0}u_0-\nabla P(\r_0)=\sqrt\r_0  g(x)
\end{equation}
are satisfied for some $g\in L^2(\mathbb{R}^2)$. Then, for any
$T>0$, there exists  a unique global classical solution
$(\r,u)(t,x)$ to the Cauchy problem \eqref{CNS}-\eqref{initial-v}
satisfying $0\le\r\leq C_2$ for some  positive constant $C_2$ and
\eqref{Jan-16-1} in Theorem \ref{theorem1}.
\end{Theorem}

\begin{Remark}
If $\l(\bar\r)<7\mu,$ one has
$\f{4(\sqrt{2+\f{\l(\bar\r)}{\mu}}-1)}{1+\f{\l(\bar\r)}{\mu}}>1.$
Then one can choose a weight $\a$ satisfying
$1<\a^2<\f{4(\sqrt{2+\f{\l(\bar\r)}{\mu}}-1)}{1+\f{\l(\bar\r)}{\mu}}$.
In this case,  the condition  $\g\leq 2\b$  in Theorems
\ref{theorem1} and \ref{theorem2} can be removed and both theorems
 hold true for any $\g>1$
and $\b>\f43$ (see \cite{JWX3} for more details).
\end{Remark}

\begin{Remark}
If $\bar\r=0,$ then
$\f{4(\sqrt{2+\f{\l(\bar\r)}{\mu}}-1)}{1+\f{\l(\bar\r)}{\mu}}=4(\sqrt2-1)$.
This is exactly same as our previous work \cite{JWX3} for the Cauchy
problem with the vanishing density at the far fields.
\end{Remark}

The rest of the paper is organized as follows. In Section 2, we
present some elementary facts  which will be used later. In Section
3, we derive a priori estimates which are needed to extend the local
solution to a global one. The sketch of proof of our main results is
given in Section 4.

 \vskip 2mm
\noindent\emph{Notations.} Throughout this paper, positive generic
constants are denoted by $c$ and $C$, which are independent of $\d$,
$m$ and $t\in[0,T]$, without confusion, and $C(\cdot)$ stands for
some generic constant(s) depending only on the quantity listed in
the parenthesis. For functional spaces, $L^{p}(\mathbb{R}^2), 1\leq
p\leq \infty$, denote the usual Lebesgue spaces on $\mathbb{R}^2$
and $\|\cdot\|_p$ denotes its $L^p$ norm. $W^{k,p}(\mathbb{R}^2)$
denotes the standard $k^{th}$ order Sobolev space and
$H^{k}(\mathbb{R}^2):=W^{k,2}(\mathbb{R}^2)$. For $1<p<\i$, the
homogenous Sobolev space $D^{k,p}(\mathbb{R}^2)$ is defined by
$D^{k,p}(\mathbb{R}^2)=\{u\in
L^1_{loc}(\mathbb{R}^2)|\|\nabla^ku\|_p<+\i\}$ with
$\|u\|_{D^{k,p}}:=\|\nabla^k u\|_p$ and
$D^k(\mathbb{R}^2):=D^{k,2}(\mathbb{R}^2)$.

\




\section{Preliminaries}

Motivated by \cite{Kazhikhov}, we introduce the following variables. First
denote the effective viscous flux by
\begin{equation}\label{flux}
F=(2\mu+\l(\r)){\rm div}u-(P(\r)-P(\bar\r)),
\end{equation}
and the vorticity by
$$
\o=\partial_{x_1}u_2-\partial_{x_2}u_1.
$$
Also, we define that
$$
H=\frac{1}{\r}(\mu\o_{x_1}+F_{x_2}),\qquad\quad
L=\frac{1}{\r}(-\mu\o_{x_2}+F_{x_1}).
$$
Then the momentum equation $\eqref{CNS}_2$ can be rewritten as
\begin{equation*}
\left\{
\begin{array}{ll}
\dot{u}_1=u_{1t}+u\cdot \nabla u_1=\frac{1}{\r}(-\mu\o_{x_2}+F_{x_1})=L,\\
\dot{u}_2=u_{2t}+u\cdot \nabla
u_2=\frac{1}{\r}(\mu\o_{x_1}+F_{x_2})=H,
\end{array}
\right.
\end{equation*}
that is,
\begin{equation*}
\dot{u}=(\dot{u}_1,\dot{u}_2)^t=(L,H)^t.
\end{equation*}
Then the effective viscous flux $F$ and the vorticity $\o$ solve the
following system:
\begin{equation*}
\left\{
\begin{array}{ll}
\o_{t}+u\cdot \nabla \o+\o{\rm div}u=H_{x_1}-L_{x_2},\\
(\f{F+P(\r)-P(\bar\r)}{2\mu+\l(\r)})_{t}+u\cdot \nabla
(\f{F+P(\r)-P(\bar\r)}{2\mu+\l(\r)})+(u_{1x_1})^2+2u_{1x_2}u_{2x_1}+(u_{2x_2})^2=H_{x_2}+L_{x_1}.
\end{array}
\right.
\end{equation*}
Due to the continuity equation $\eqref{CNS}_1$, it holds that
\begin{equation}\label{F-omega}
\left\{
\begin{array}{ll}
\o_{t}+u\cdot \nabla \o+\o{\rm div}u=H_{x_1}-L_{x_2},\\
F_{t}+u\cdot \nabla
F-\r(2\mu+\l(\r))[F(\f{1}{2\mu+\l(\r)})^\prime+(\f{P(\r)-P(\bar\r)}{2\mu+\l(\r)})^\prime]{\rm
div}u\\
\qquad+(2\mu+\l(\r))[(u_{1x_1})^2+2u_{1x_2}u_{2x_1}+(u_{2x_2})^2]=(2\mu+\l(\r))(H_{x_2}+L_{x_1}).
\end{array}
\right.
\end{equation}
Furthermore,  the system for $(H,L)$ can be derived as
\begin{equation*}
\left\{
\begin{array}{ll}
\r H_{t}+\r u\cdot \nabla H-\r H{\rm div}u+u_{x_2}\cdot\nabla F+\mu
u_{x_1}\cdot\nabla\o+\mu(\o{\rm
div}u)_{x_1}\\
\qquad
-\big\{\r(2\mu+\l(\r))[F(\f{1}{2\mu+\l(\r)})^\prime+(\f{P(\r)-P(\bar\r)}{2\mu+\l(\r)})^\prime]{\rm
div}u\big\}_{x_2}\\
\qquad+\big\{(2\mu+\l(\r))[(u_{1x_1})^2+2u_{1x_2}u_{2x_1}+(u_{2x_2})^2]\big\}_{x_2}\\
\qquad=[(2\mu+\l(\r))(H_{x_2}+L_{x_1})]_{x_2}+\mu(H_{x_1}-L_{x_2})_{x_1},\\
\r L_{t}+\r u\cdot \nabla L-\r L{\rm div}u+u_{x_1}\cdot\nabla F-\mu
u_{x_2}\cdot\nabla\o-\mu(\o{\rm
div}u)_{x_2}\\
\qquad
-\big\{\r(2\mu+\l(\r))[F(\f{1}{2\mu+\l(\r)})^\prime+(\f{P(\r)-P(\bar\r)}{2\mu+\l(\r)})^\prime]{\rm
div}u\big\}_{x_1}\\
\qquad+\big\{(2\mu+\l(\r))[(u_{1x_1})^2+2u_{1x_2}u_{2x_1}+(u_{2x_2})^2]\big\}_{x_1}\\
\qquad=[(2\mu+\l(\r))(H_{x_2}+L_{x_1})]_{x_1}-\mu(H_{x_1}-L_{x_2})_{x_2}.
\end{array}
\right.
\end{equation*}
In the following, we will utilize the above systems in different
steps. Note that these systems are equivalent to each other for the
smooth solution to the original system \eqref{CNS}.

Several elementary Lemmas are needed later. The first one is the
various Gagliardo-Nirenberg inequalities.

\begin{Lemma}\label{lemma1}

\begin{itemize}
\item[(1)]
$\forall h\in W^{1,m}(\mathbb{R}^2)\cap L^r(\mathbb{R}^2)$, it holds
that
\begin{equation*}
\|h\|_q\leq C\|\nabla h\|_m^\t\|h\|_r^{1-\t},
\end{equation*}
where $\t=(\f1r-\f1q)(\f1r-\f1m+\f12)^{-1}$, and if $m<2,$ then $q$
is between $r$ and $\f{2m}{2-m}$, that is, $q\in[r,\f{2m}{2-m}]$ if
$r<\f{2m}{2-m}$, $q\in[\f{2m}{2-m},r]$ if $r\geq\f{2m}{2-m},$ if
$m=2,$ then $q\in[r,+\i)$, if $m>2$, then $q\in[r,+\i].$
\item[(2)](Best constant for the Gagliardo-Nirenberg inequality)

$\forall h\in \mathbb{D}^m(\mathbb{R}^2)\doteq\Big\{h\in
L^{m+1}(\mathbb{R}^2)\Big|\nabla h\in L^2(\mathbb{R}^2),h\in
L^{2m}(\mathbb{R}^2)\Big\}$ with $m>1$, it holds that
\begin{equation*}
\|h\|_{2m}\leq A_m\|\nabla h\|_2^\t\|u\|_{m+1}^{1-\t},
\end{equation*}
where $\t=\f12-\f{1}{2m}$ and
$$
A_m=\Big(\f{m+1}{2\pi}\Big)^{\f{\t}{2}}\Big(\f{2}{m+1}\Big)^{\f{1}{2m}}\leq
C m^{\f14}
$$
with the positive constant $C$ independent of $m$, and $A_m$ is the
optimal constant.

\item[(3)] $\forall h\in W^{1,m}(\mathbb{R}^2)$ with
$1\leq m<2,$ then
\begin{equation*}
\|h\|_{\f{2m}{2-m}}\leq C(2-m)^{-\f12}\|\nabla h\|_m,
\end{equation*}
where the positive constant $C$ is independent of $m.$

\end{itemize}
\end{Lemma}

{\bf Proof:} The proof of (1) can be found in \cite{Kazhikhov} while the
proof of (2) can be found in \cite{Del-Pino}.
The proof of (3) can be found in \cite{GT}.

The following Lemma is  the Caffarelli-Kokn-Nirenberg weighted
inequalities, which is crucial to the weighted estimates in the two-dimensional
Cauchy problem.

\begin{Lemma}\label{lemma2}
\begin{itemize}
\item [(1)] $\forall h\in C^\i_0(\mathbb{R}^2)$, it holds that
\begin{equation*}
\||x|^\kappa h\|_r\leq C\||x|^\a |\nabla h|\|^\t_p~\||x|^\b
h\|^{1-\t}_q
\end{equation*}
where $1\leq p,q<\i, 0<r<\i, 0\leq \t\leq 1,
\f1p+\f{\a}{2}>0,\f1q+\f\b2>0,\f1r+\f\kappa2>0$ and satisfying
\begin{equation*}
\f1r+\f\kappa2=\t(\f1p+\f{\a-1}{2})+(1-\t)(\f1q+\f\b2),
\end{equation*}
and
$$
\kappa=\t\sigma+(1-\t)\b,
$$
with $0\leq \a-\sigma$ if $\t>0$ and $0\leq \a-\sigma\leq 1$ if
$\t>0$ and $\f1p+\f{\a-1}{2}=\f1r+\f\kappa2.$

\item [(2)](Best constant for Caffarelli-Kohn-Nirenberg inequality)

$\forall h\in C^\i_0(\mathbb{R}^2)$, it holds that
\begin{equation}\label{b-CKN}
\||x|^b h\|_p\leq C_{a,b}\||x|^a\nabla h\|_2
\end{equation}
where $a>0, a-1\leq b\leq a$ and $p=\f{2}{a-b}$.  If $b=a-1$, then
$p=2$ and the best constant in the inequality \eqref{b-CKN} is
$$C_{a,b}=C_{a,a-1}=a.$$
\end{itemize}
\end{Lemma}
{\bf Proof:} The proof of (1) can be found in \cite{CKN} while the
proof of (2) can be found in \cite{CW}.

The following lemma is known and the proof is referred to \cite{JWX3}.
\begin{Lemma}\label{lemma3}
\begin{itemize}
\item [(1)] It holds that for $1<p<\i$ and $u\in C_0^\i(\mathbb{R}^2)$,
\begin{equation*}
\|\nabla u\|_p\leq C(\|{\rm div} u\|_p+\|\o\|_p);
\end{equation*}

\item [(2)] It holds that for $1<p<\i$, $-2<\a<2(p-1)$ and $u\in C_0^\i(\mathbb{R}^2)$,
\begin{equation*}
\||x|^{\f\a p}|\nabla u|\|_p\leq C(\||x|^{\f\a p}{\rm div}
u\|_p+\||x|^{\f\a p}\o\|_p).
\end{equation*}
\end{itemize}
\end{Lemma}

\section{A priori estimates}
\setcounter{equation}{0}

In this section, we will obtain various a priori estimates and a upper bound of the density.

\underline{Step 1. Elementary energy estimates:}
\begin{Lemma}\label{lemma-ee}
There exists a positive constant $C$ depending on $(\r_0,u_0)$, such
that
\begin{equation*}
\sup_{t\in[0,T]}\big(\|\sqrt\r
u\|^2_2+\|\Psi(\r,\bar\r)\|_1\big)+\int_0^T\big(\|\nabla
u\|_2^2+\|\o\|_2^2+\|(2\mu+\l(\r))^\f12{\rm div} u\|_2^2\big)dt\leq
C.
\end{equation*}
\end{Lemma}
{\bf Proof:} Multiplying the equation $\eqref{CNS}_2$ by $u$, the continuity equation $\eqref{CNS}_1$ by
$\f{\g}{\g-1}\r^{\g-1}$, then summing the resulting equations,  and using the
continuity equation $\eqref{CNS}_1$,  yield that
\begin{equation}\label{e1}
\begin{array}{ll}
\di \big[\r\f{|u|^2}{2}+\Psi(\r,\bar\r)\big]_t+ {\rm div}\big[\r u\f{|u|^2}{2}+\Psi(\r,\bar\r)u+(P(\r)-P(\bar\r))u\big]\\
\di ={\rm div}\big[\mu\nabla\f{|u|^2}{2}+(\mu+\l(\r))({\rm div}
u)u\big]-\mu|\nabla u|^2-(\mu+\l(\r))({\rm div} u)^2.
\end{array}
\end{equation}
Therefore, integrating the above equality over
$[0,t]\times\mathbb{R}^2$ with respect to $t$ and $x$ and noting
that
\begin{equation*}
\int \big[\mu|\nabla u|^2+(\mu+\l(\r))({\rm div} u)^2\big]dx =\int \big[\mu\o^2+(2\mu+\l(\r))({\rm div} u)^2\big]dx,
\end{equation*}
complete the proof of Lemma \ref{lemma-ee}. $\hfill\Box$

\underline{Step 2. Weighted energy estimates:}

The following weighted energy estimates are fundamental and crucial
in our analysis.
\begin{Lemma}\label{lemma-wee}
For $\a>0$ satisfying
$\a^2<\f{4(\sqrt{2+\f{\l(\bar\r)}{\mu}}-1)}{1+\f{\l(\bar\r)}{\mu}}$ and $\g\leq 2\b$,
it holds that for sufficiently large $m>1$ and $\forall t\in[0,T]$,
\begin{equation}\label{We-l1}
\begin{array}{ll}
\di \int_{\mathbb{R}^2}|x|^\a\big[\r|u|^2+\Psi(\r,\bar\r)\big](t,x)dx+\int_0^t\big[\||x|^{\f\a2}\nabla u\|_2^2(s)+\||x|^{\f\a2}{\rm div} u\|_2^2(s)+\||x|^{\f\a2}\sqrt{\l(\r)}{\rm div} u\|_2^2(s)\big]ds\\
\di \leq
C_\a\Big[1+\int_0^t(\|\r-\bar\r\|^\b_{2m\b+1}(s)+1)(\|\nabla
u\|_2^2(s)+1)ds\Big],
\end{array}
\end{equation}
where the positive constant $C_\a$ may depend on $\a$ but is
independent of $m$.
\end{Lemma}
{\bf Proof:} Multiplying the equality \eqref{e1} by $|x|^\a$ yields
that
\begin{equation}\label{e2}
\begin{array}{ll}
\di \big[|x|^\a(\r\f{|u|^2}{2}+\Psi(\r,\bar\r))\big]_t+\big[\mu|\nabla u|^2+(\mu+\l(\r))({\rm div} u)^2\big]|x|^\a\\
\di =-{\rm div}\big[|x|^\a\big(\r u\f{|u|^2}{2}+\Psi(\r,\bar\r)u+(P(\r)-P(\bar\r))u\big)\big]+{\rm div}\big[\big(\mu\nabla\f{|u|^2}{2}+(\mu+\l(\r))({\rm div} u)u\big)|x|^\a\big]\\
\di \quad+\big[\r
u\f{|u|^2}{2}+\Psi(\r,\bar\r)u+(P(\r)-P(\bar\r))\big]u\cdot\nabla(|x|^\a)-\big[\mu\nabla\f{|u|^2}{2}+(\mu+\l(\r))({\rm
div} u)u\big]\cdot\nabla(|x|^\a).
\end{array}
\end{equation}
Integrating the above equation \eqref{e2} with respect to $x$ over
$\mathbb{R}^2$ yields that
\begin{equation}\label{e3}
\begin{array}{ll}
\di \f{d}{dt}\int |x|^\a\big[\r\f{|u|^2}{2}+\Psi(\r,\bar\r)\big](t,x)dx+\Big[\mu\||x|^{\f\a2}\nabla u\|_2^2
+\mu\||x|^{\f\a2}{\rm div} u\|_2^2+\||x|^{\f\a2}\sqrt{\l(\r)}{\rm div} u\|_2^2\Big](t) \\
\di =\int\big[\r
u\f{|u|^2}{2}+\Psi(\r,\bar\r)u+(P(\r)-P(\bar\r))\big]u\cdot\nabla(|x|^\a)dx\\
\di\qquad\qquad\qquad  -\int
\big[\mu\nabla\f{|u|^2}{2}+(\mu+\l(\r))({\rm div}
u)u\big]\cdot\nabla(|x|^\a)dx.
\end{array}
\end{equation}

Now we estimate the terms on the right hand side of \eqref{e3}.
First, it holds that
\begin{equation}\label{WE1}
\begin{array}{ll}
\di |\int\r \f{|u|^2}{2}u\cdot\nabla(|x|^\a) dx|\di =|\int\f{|u|^2}{2}\big((\sqrt\r-\sqrt{\bar\r})+\sqrt{\bar\r}\big)\sqrt\r u\cdot \nabla(|x|^\a) dx|\\
\di\quad\leq |\int\f{|u|^2}{2}\big(\sqrt\r-\sqrt{\bar\r}\big)\sqrt\r
u\cdot \nabla(|x|^\a) dx|+\sqrt{\bar\r}~|\int\f{|u|^2}{2}\sqrt\r
u\cdot \nabla(|x|^\a) dx|:=I_{11}+I_{12}.
\end{array}
\end{equation}
Then, it follows that
\begin{equation}\label{WE110}
\begin{array}{ll}
\di I_{11}=|\int\f{|u|^2}{2}\big(\sqrt\r-\sqrt{\bar\r}\big)\big({\bf
1}|_{\{0\leq\r\leq2\bar\r\}}+{\bf 1}|_{\{\r>2\bar\r\}}\big)\sqrt\r
u\cdot \nabla(|x|^\a) dx|\\
\di \leq C\|\sqrt\r u\|_2\Big[\|(\sqrt\r-\sqrt{\bar\r}){\bf
1}|_{\{0\leq\r\leq2\bar\r\}}\|_{p_1}\||x|^{\a-1}|u|^2\|_{q_1}+\|(\sqrt\r-\sqrt{\bar\r}){\bf 1}|_{\{\r>2\bar\r\}}\|_{2\g}\||x|^{\a-1}|u|^2\|_{\f{2\g}{\g-1}}\Big]\\
\di \leq C\big[\|(\r-\bar\r){\bf
1}|_{\{0\leq\r\leq2\bar\r\}}\|_{p_1}\||x|^{\f{\a-1}{2}}u\|^2_{2q_1}+\|\Psi(\r,\bar\r){\bf 1}|_{\{\r>2\bar\r\}}\|^\f1{2\g}_{1}\||x|^{\f{\a-1}{2}}u\|^2_{\f{4\g}{\g-1}}\big]\\
\di  \leq C\big[\|(\r-\bar\r){\bf
1}|_{\{0\leq\r\leq2\bar\r\}}\|_{p_1}\|\nabla u\|_2^{2\t_1}\||x|^{\f{\a}{2}}\nabla u\|^{2(1-\t_1)}_{2}+\|\nabla u\|_2^{\f{2}{\a\g}}\||x|^{\f{\a}{2}}\nabla u\|^{2(1-\f{1}{\a\g})}_{2}\big]\\
\di \leq \sigma\||x|^{\f{\a}{2}}\nabla
u\|^{2}_{2}+C_\sigma\big[1+\|(\r-\bar\r){\bf
1}|_{\{0\leq\r\leq2\bar\r\}}\|^{\f1{\t_1}}_{p_1}\big]\|\nabla u\|_2^2,
\end{array}
\end{equation}
where and in the sequel $\sigma >0$ is a small constant to be
determined, $C_\sigma$ is a positive constant depending on $\sigma$.
By the H${\rm\ddot{o}}$lder inequality and the
Caffarelli-Kohn-Nirenberg inequality in  Lemma \ref{lemma2} (1), the
positive constants $p_1>2, q_1>2, \t_1\in(0,1]$ in the above
inequality  \eqref{WE11} satisfying
$$
\f{1}{p_1}+\f{1}{q_1}=\f12,
$$
and
$$
\f{1}{2q_1}+\f{\f{\a-1}{2}}{2}=\t_1(\f12+\f{0-1}{2})+(1-\t_1)(\f{1}{2}+\f{\f{\a}{2}-1}{2})=\f{\a}{4}(1-\t_1).
$$
The combination of the above two equalities yields that
\begin{equation}\label{W-p1}
p_1=\f{2}{\a\t_1},
\end{equation}
with $\a>0$, $\t_1\in(0,1)$ and $p_1>2$. Therefore, it holds that
$$
\|(\r-\bar\r){\bf
1}|_{\{0\leq\r\leq2\bar\r\}}\|^{\f1{\t_1}}_{p_1}\leq C\|(\r-\bar\r){\bf
1}|_{\{0\leq\r\leq2\bar\r\}}\|^{\f1{\t_1}}_2\leq C \|\Psi(\r,\bar\r){\bf
1}|_{\{0\leq\r\leq2\bar\r\}}\|^{\f1{2\t_1}}_1\leq C,
$$
which together with \eqref{WE110} gives that
\begin{equation}\label{WE11}
\di I_{11} \leq \sigma\||x|^{\f{\a}{2}}\nabla
u\|^{2}_{2}+C_\sigma\|\nabla u\|_2^2.
\end{equation}
Then, one can obtain
\begin{equation}\label{WE12}
\begin{array}{ll}
\di\quad I_{12}\leq \f{\a}{2}\sqrt{\bar\r}\|\sqrt\r u|x|^{\f\a2}\|_2\||x|^{\f{\f\a2-1}{2}}|u|\|^2_{4}\\
\di \qquad\leq C\|\sqrt\r u|x|^{\f\a2}\|_2\|\nabla u\|_2\||x|^{\f{\a}{2}}\nabla u\|_{2}\\
\di \qquad \leq \sigma\||x|^{\f{\a}{2}}\nabla
u\|^{2}_{2}+C_\sigma\|\sqrt\r u|x|^{\f\a2}\|_2^2\|\nabla u\|_2^2.
\end{array}
\end{equation}
Then it holds that
\begin{equation}\label{WE2}
\begin{array}{ll}
\di |\int\big[\Psi(\r,\bar\r)+(P(\r)-P(\bar\r))\big]u\cdot\nabla(|x|^\a) dx|\\
\di\quad
=|\int\big[\Psi(\r,\bar\r)+(P(\r)-P(\bar\r))\big]\big({\bf
1}|_{\{0\leq\r\leq2\bar\r\}}+{\bf 1}|_{\{\r>2\bar\r\}}\big)u\cdot\nabla(|x|^\a) dx|\\
\di\quad \leq C|\int\big[|\r-\bar\r|{\bf
1}|_{\{0\leq\r\leq2\bar\r\}}+|\r-\bar\r|^\g{\bf 1}|_{\{\r>2\bar\r\}}\big]u\cdot\nabla(|x|^\a) dx|\\
\di \quad\leq C\big[\|(\r-\bar\r){\bf
1}|_{\{0\leq\r\leq2\bar\r\}}\|_{2}\||x|^{\a-1}u\|_2+\|(\r-\bar\r){\bf 1}|_{\{\r>2\bar\r\}}\|^\g_{\g
p_2}\||x|^{\a-1}u\|_{q_2}\big]\\
\di \quad \leq C\big[\|\Psi(\r,\bar\r){\bf
1}|_{\{0\leq\r\leq2\bar\r\}}\|^{\f12}_{1}\|\nabla u\|_2^{\f12}\|\nabla u |x|^{\f\a2}\|_2^{\f12}+\|(\r-\bar\r){\bf 1}|_{\{\r>2\bar\r\}}\|^\g_{\g
p_2}\|\nabla u\|_2^{\t_2}\||x|^{\f{\a}{2}}\nabla u\|^{1-\t_2}_{2}\big]\\
\di \quad \leq \sigma\||x|^{\f{\a}{2}}\nabla u\|^{2}_{2}+C_\sigma
\big[1+\|(\r-\bar\r){\bf 1}|_{\{\r>2\bar\r\}}\|^{\f{2\g}{1+\t_2}}_{\g
p_2}\big](\|\nabla u\|_2^2+1),
\end{array}
\end{equation}
By the H${\rm\ddot{o}}$lder inequality and the
Caffarelli-Kohn-Nirenberg inequality in  Lemma \ref{lemma2} (1), the
positive constants $p_2>1, q_2>1, \t_2\in(0,1]$ in the above
inequality \eqref{WE2} satisfying
$$
\f{1}{p_2}+\f{1}{q_2}=1,
$$
$$
\f{1}{q_2}+\f{\a-1}{2}=\t_2(\f12+\f{0-1}{2})+(1-\t_2)(\f{1}{2}+\f{\f{\a}{2}-1}{2})=\f{\a}{4}(1-\t_2).
$$
The combination of the above three equalities yields that
\begin{equation}\label{W-p2}
p_2=\f{4}{2+\a(1+\t_2)},
\end{equation}
with the parameters $\a>0$, $\t_2\in(0,1)$
and  $p_2>1$. Note that $p_2>1$ is equivalent to the condition that
\begin{equation}\label{c11}
\f{\a}{2}(1+\t_2)<1.
\end{equation}
Then one can compute that
\begin{equation}\label{WE3}
\begin{array}{ll}
\di |-\int \mu\nabla\f{|u|^2}{2}\cdot\nabla(|x|^\a)dx|=\mu\a|\int u\cdot\nabla u\cdot x|x|^{\a-2}dx|\\
\di\quad \leq \mu\a\||x|^{\f\a2}\nabla u\|_2\||x|^{\f\a2-1}
u\|_2\leq \f{\mu\a^2}{2}\||x|^{\f\a2}\nabla u\|_2^2,
\end{array}
\end{equation}
where in the last inequality one has used the best constant $\f\a2$
for the Caffarelli-Kohn-Nirenberg inequality in Lemma \ref{lemma2}
(2). Similarly, it holds that
\begin{equation}\label{WE4}
\begin{array}{ll}
\di |-\int \mu({\rm div}u)u\cdot\nabla(|x|^\a)dx|=\mu\a|\int ({\rm div}u) |x|^{\a-2}u\cdot x dx|\\
\di\quad \leq \mu\a\||x|^{\f\a2}{\rm div}u\|_2\||x|^{\f\a2-1}
u\|_2\leq \f{\mu\a^2}{2}\||x|^{\f\a2}{\rm
div}u\|_2\||x|^{\f\a2}\nabla u\|_2.
\end{array}
\end{equation}
Then it follows that
\begin{equation}\label{WE5}
\begin{array}{ll}
\di |-\int \l(\r)({\rm div}u)u\cdot\nabla(|x|^\a)dx|
\di =\a|\int \sqrt{\l(\r)}({\rm div}u)\big[(\sqrt{\l(\r)}-\sqrt{\l(\bar\r)})+\sqrt{\l(\bar\r)}\big] |x|^{\a-2}u\cdot x dx|\\
\di \leq \a|\int \sqrt{\l(\r)}({\rm
div}u)\big(\sqrt{\l(\r)}-\sqrt{\l(\bar\r)}\big) |x|^{\a-2}u\cdot x
dx|+\sqrt{\l(\bar\r)}\a |\int \sqrt{\l(\r)}({\rm div}u)
|x|^{\a-2}u\cdot x dx|\\
\di:=I_{21}+I_{22}.
\end{array}
\end{equation}
It holds that
\begin{equation}\label{WE510}
\begin{array}{ll}
\di I_{21}\leq \a\|\sqrt{\l(\r)}|x|^{\f\a2}{\rm div}u\|_2\|\sqrt{\l(\r)}-\sqrt{\l(\bar\r)}\|_{p_3}\||x|^{\f\a2-1} u\|_{q_3}\\
\di \leq C\|\sqrt{\l(\r)}|x|^{\f\a2}{\rm div}u\|_2\big[\|(\r-\bar\r){\bf
1}|_{\{0\leq\r\leq2\bar\r\}}\|_{p_3}+\|(\r-\bar\r){\bf 1}|_{\{\r>2\bar\r\}}\|^{\f\b2}_{\f{\b p_3}{2}}\big]\|\nabla u\|_2^{\t_3}\||x|^{\f\a2} \nabla u\|_{2}^{1-\t_3}\\
\di \leq \sigma\big[\|\sqrt{\l(\r)}|x|^{\f\a2}{\rm
div}u\|_2^2+\||x|^{\f\a2} \nabla
u\|_{2}^{2}\big]\\
\di\qquad\qquad\qquad +C_\sigma\big[\|(\r-\bar\r){\bf
1}|_{\{0\leq\r\leq2\bar\r\}}\|^{\f2{\t_3}}_{p_3}+\|(\r-\bar\r){\bf 1}|_{\{\r>2\bar\r\}}\|^{\f\b{\t_3}}_{\f{\b
p_3}{2}}\big]\|\nabla u\|_2^{2}.
\end{array}
\end{equation}
By the H${\rm\ddot{o}}$lder inequality and the
Caffarelli-Kohn-Nirenberg inequality in  Lemma \ref{lemma2} (1), the
positive constants $p_3>2, q_3>2, \t_3\in(0,1]$ in the above
inequality \eqref{WE51} satisfying
$$
\f{1}{p_3}+\f{1}{q_3}=\f12,
$$
$$
\f{1}{q_3}+\f{\f{\a}{2}-1}{2}=\t_3(\f12+\f{0-1}{2})+(1-\t_3)(\f{1}{2}+\f{\f{\a}{2}-1}{2})=\f{\a}{4}(1-\t_3).
$$
The combination of the above three equalities yields that
\begin{equation}\label{W-p3}
p_3=\f{4}{\a\t_3}.
\end{equation}
with $\a>0$, $\t_3\in(0,1)$ and  $p_3>2$. Therefore, it holds that
$$
\|(\r-\bar\r){\bf
1}|_{\{0\leq\r\leq2\bar\r\}}\|^{\f2{\t_3}}_{p_3}\leq C\|(\r-\bar\r){\bf
1}|_{\{0\leq\r\leq2\bar\r\}}\|^{\f2{\t_3}}_2\leq C \|\Psi(\r,\bar\r){\bf
1}|_{\{0\leq\r\leq2\bar\r\}}\|^{\f1{\t_3}}_1\leq C,
$$
which together with \eqref{WE510} gives that
\begin{equation}\label{WE51}
I_{21}\leq \sigma\big[\|\sqrt{\l(\r)}|x|^{\f\a2}{\rm
div}u\|_2^2+\||x|^{\f\a2} \nabla
u\|_{2}^{2}\big]+C_\sigma\big[1+\|(\r-\bar\r){\bf 1}|_{\{\r>2\bar\r\}}\|^{\f\b{\t_3}}_{\f{\b
p_3}{2}}\big]\|\nabla u\|_2^{2}.
\end{equation}
Meanwhile, it holds that
\begin{equation}\label{WE52}
\begin{array}{ll}
\di\quad I_{22}\leq \sqrt{\l(\bar\r)}\a\|\sqrt{\l(\r)}|x|^{\f\a2}{\rm div}u\|_2\||x|^{\f\a2-1} u\|_{2}\\
\di \qquad \leq
\f{\a^2\sqrt{\l(\bar\r)}}{2}\|\sqrt{\l(\r)}|x|^{\f\a2}{\rm
div}u\|_2\||x|^{\f\a2}\nabla u\|_{2}.
\end{array}
\end{equation}
Substituting \eqref{WE11} and \eqref{WE12} into \eqref{WE1},
\eqref{WE51} and \eqref{WE52} into \eqref{WE5} and then substituting
the resulting \eqref{WE1}, \eqref{WE5} and  \eqref{WE2},
 \eqref{WE3} and \eqref{WE4}  into \eqref{e3} yield
that
\begin{equation}\label{e4}
\begin{array}{ll}
\di \f{d}{dt}\int |x|^\a\big[\r\f{|u|^2}{2}+\Psi(\r,\bar\r)\big](t,x)dx+J(t)\leq \sigma\big[\|\sqrt{\l(\r)}|x|^{\f\a2}{\rm
div}u\|_2^2+4\||x|^{\f\a2} \nabla
u\|_{2}^{2}\big]\\
\di ~+C_\sigma\big[1+\|\sqrt\r
u|x|^{\f\a2}\|_2^2+\|(\r-\bar\r){\bf 1}|_{\{\r>2\bar\r\}}\|^{\f{2\g}{1+\t_2}}_{p_2\g}+\|(\r-\bar\r){\bf 1}|_{\{\r>2\bar\r\}}\|^{\f\b{\t_3}}_{\f{\b
p_3}{2}}\big](\|\nabla u\|_2^{2}+1),
\end{array}
\end{equation}
where $\t_i\in(0,1]$, $p_i~(i=1,2,3)$ are given in \eqref{W-p1},
\eqref{W-p2} and \eqref{W-p3}, respectively, and $p_1, p_3>2$ and
$p_2>1$ and
\begin{equation}
\label{J}
\begin{array}{ll}
\di J(t)=\mu(1-\f{\a^2}{2})\||x|^{\f\a2}\nabla u\|_2^2-\f{\mu\a^2}{2}\||x|^{\f\a2}\nabla u\|_2\||x|^{\f\a2}{\rm div} u\|_2+\mu\||x|^{\f\a2}{\rm div} u\|_2^2\\
\di\qquad +\||x|^{\f\a2}\sqrt{\l(\r)}{\rm div} u\|_2^2(t)
-\f{\a^2\sqrt{\l(\bar\r)}}{2}\||x|^{\f\a2}\sqrt{\l(\r)}{\rm div}
u\|_2\||x|^{\f\a2}\nabla
u\|_2.
\end{array}
\end{equation}
The corresponding matrix of the above quadratic term \eqref{J} is
$$
A=\left(
\begin{array}{ccc}
\di \mu(1-\f{\a^2}{2})&\di -\f{\mu\a^2}{4}&\di-\f{\a^2\sqrt{\l(\bar\r)}}{4}\\[2mm]
\di -\f{\mu\a^2}{4}&\mu&0\\[2mm]
\di -\f{\a^2\sqrt{\l(\bar\r)}}{4}&0&1
\end{array}
\right).
$$
The matrix $A$  is positively definite if and only if all the principal minor determinant of $A$ is positive, that is,
$$
\begin{array}{ll}
\di \mu(1-\f{\a^2}{2})>0,\qquad
\left|
\begin{array}{ccc}
\di \mu(1-\f{\a^2}{2})&\di -\f{\mu\a^2}{4}\\[2mm]
\di -\f{\mu\a^2}{4}&\mu\\
\end{array}
\right|>0,~~{\rm and}~~~
\di \left|
\begin{array}{ccc}
\di \mu(1-\f{\a^2}{2})&\di -\f{\mu\a^2}{4}&\di-\f{\a^2\sqrt{\l(\bar\r)}}{4}\\[2mm]
\di -\f{\mu\a^2}{4}&\mu&0\\[2mm]
\di -\f{\a^2\sqrt{\l(\bar\r)}}{4}&0&1
\end{array}
\right|>0.
\end{array}
$$
Therefore, if the weight $\a$ satisfies
\begin{equation}\label{alpha}
0<\a^2<\f{4(\sqrt{2+\f{\l(\bar\r)}{\mu}}-1)}{1+\f{\l(\bar\r)}{\mu}},
\end{equation}
then the matrix $A$ is positively definite, and then there exists a
positive constant $C_\a$ such that
\begin{equation}\label{quad}
J(t)\geq C^{-1}_\a\Big[\||x|^{\f\a2}\nabla
u\|_2^2(t)+\||x|^{\f\a2}{\rm div}
u\|_2^2(t)+\||x|^{\f\a2}\sqrt{\l(\r)}{\rm div} u\|_2^2(t)\Big].
\end{equation}
Consequently, if the weight $\a$ satisfies \eqref{alpha}, then
substituting \eqref{quad} into \eqref{e4} and choosing $\sigma$
suitably small yield that
\begin{equation}\label{e5-1}
\begin{array}{ll}
\di \f{d}{dt}\int |x|^\a\big[\r\f{|u|^2}{2}+\Psi(\r,\bar\r)\big](t,x)dx+\f{\mu}{2}C^{-1}_\a\Big[\||x|^{\f\a2}\nabla u\|_2^2+\||x|^{\f\a2}{\rm div} u\|_2^2+\||x|^{\f\a2}\sqrt{\l(\r)}{\rm div} u\|_2^2(t)\Big]\\
\di \leq C\big[1+\|\sqrt\r
u|x|^{\f\a2}\|_2^2+
+\|(\r-\bar\r){\bf 1}|_{\{\r>2\bar\r\}}\|^{\f{2\g}{1+\t_2}}_{p_2\g}+\|(\r-\bar\r){\bf 1}|_{\{\r>2\bar\r\}}\|^{\f\b{\t_3}}_{\f{\b
p_3}{2}}\big](\|\nabla u\|_2^{2}+1).
\end{array}
\end{equation}
Now choose $m>1$ sufficiently large such that
$$2m\b+1\geq \max\Big\{p_2\g, \f{\b p_3}{2}\Big\}.$$
Then, it holds that
\begin{equation}\label{Im2}
\|(\r-\bar\r){\bf 1}|_{\{\r>2\bar\r\}}\|^{\f{2\g}{1+\t_2}}_{p_2\g}\leq \|(\r-\bar\r){\bf 1}|_{\{\r>2\bar\r\}}\|^{\f{2a_2\g}{1+\t_2}}_{\g}\|\r-\bar\r\|^{\f{2(1-a_2)\g}{1+\t_2}}_{2m\b+1}\leq C\|\Psi(\r,\bar\r)\|^{\f{2a_2}{1+\t_2}}_{1}\|\r-\bar\r\|^{\f{2\g(1-a_2)}{1+\t_2}}_{2m\b+1},
\end{equation}
with $a_2\in(0,1)$ satisfying
$$
\f{a_2}{\g}+\f{1-a_2}{2m\b+1}=\f{1}{p_2\g}=\f{2+\a(1+\t_2)}{4\g},
$$
which implies that
$$
a_2=\f{(2+\a(1+\t_2))(2m\b+1)-4\g}{4(2m\b+1-\g)}\rightarrow \f{2+\a(1+\t_2)}{4},~~{\rm as}~m\rightarrow+\i.
$$
The following restriction should be imposed to \eqref{Im2}
$$
\f{2\g(1-a_2)}{1+\t_2}\leq \b,
$$
which is satisfied provided
\begin{equation}\label{a2}
(1+\t_2)(\f{\b}{\g}+\f{\a}{2})>1
\end{equation}
and $m\gg1.$
Since $\g\leq 2\b$, then $\f\b\g\geq\f12,$ thus
 one can choose $\t_2\in(0,1)$ such that $(1+\t_2)\f{\b}{\g}\geq1$ and $\f{\a}{2}(1+\t_2)<1$, and thus satisfies
the restrictions \eqref{a2} and \eqref{c11} if $m\gg1$. If $\l(\bar\r)<7\mu,$ then $\f{4(\sqrt{2+\f{\l(\bar\r)}{\mu}}-1)}{1+\f{\l(\bar\r)}{\mu}}>1.$ Thus we can choose the weight $\a>0$ satisfying $1<\a^2<\f{4(\sqrt{2+\f{\l(\bar\r)}{\mu}}-1)}{1+\f{\l(\bar\r)}{\mu}}$. In this case, one can choose $\t_2\in(0,1)$ satisfying
the restrictions \eqref{c11} and \eqref{a2} for any fixed $\g,\b>1$, that is, the condition  $\g\leq 2\b$ in the Theorem \ref{theorem2} can be removed as in Remark 2. Then it follows from
\eqref{Im2} that
\begin{equation}\label{Im2-1}
\|(\r-\bar\r){\bf 1}|_{\{\r>2\bar\r\}}\|_{p_2\g}^{\f{2\g}{1+\t_2}}\leq C(\|\r-\bar\r\|^{\b}_{2m\b+1}+1)
\end{equation}
with the positive constant $C$ independent of $m$.

Similarly, one has
\begin{equation}\label{Im3-1}
\begin{array}{ll}
\di \|(\r-\bar\r){\bf 1}|_{\{\r>2\bar\r\}}\|^{\f{\b}{\t_3}}_{\f{\b
p_3}{2}}\leq \|(\r-\bar\r){\bf 1}|_{\{\r>2\bar\r\}}\|^{\f{\b}{\t_3}a_3}_{1}\|\r-\bar\r\|^{\f{\b}{\t_3}(1-a_3)}_{2m\b+1}\\
 \di \leq \|\Psi(\r,\bar\r){\bf 1}|_{\{\r>2\bar\r\}}\|^{\f{\b}{\t_3}a_3}_{1}\|\r-\bar\r\|^{\f{\b}{\t_3}(1-a_3)}_{2m\b+1},
\end{array}
\end{equation}
with with $a_3\in(0,1)$ satisfying
$$
\f{a_3}{1}+\f{1-a_3}{2m\b+1}=\f{2}{p_3\b}=\f{\a\t_3}{2\b},
$$
which implies that
$$
a_3=\f{\a\t_3(2m\b+1)-2\b}{4m\b^2}\rightarrow \f{\a\t_3}{2\b},~~{\rm as}~m\rightarrow+\i.
$$
The following restriction should be imposed to \eqref{Im3-1}
$$
\f{1-a_3}{\t_3}\leq 1,
$$
which is satisfied provided we choose $m\gg1$ and $\t_3\in(0,1)$ such that
\begin{equation*}
\t_3(1+\f{\a}{2\b})>1.
\end{equation*}
Then it follows from
\eqref{Im3-1} that
\begin{equation}\label{Im3-2}
\|(\r-\bar\r){\bf 1}|_{\{\r>2\bar\r\}}\|^{\f{\b}{\t_3}}_{\f{\b
p_3}{2}}\leq C(\|\r-\bar\r\|^{\b}_{2m\b+1}+1)
\end{equation}
with the positive constant $C$ independent of $m$.
Substituting \eqref{Im2-1} and \eqref{Im3-2} into
\eqref{e5-1}, then integrating the resulting inequality over $[0,t]$
with $t\in[0,T]$ and using Gronwall inequality yield the estimate \eqref{We-l1} in Lemma
\ref{lemma-wee}. $\hfill\Box$

\underline{Step 3. Density estimates:}

Applying the operator $div$ to the momentum equation
$\eqref{CNS}_2$, it holds that
\begin{equation}\label{vis-f}
[{\rm div}(\r u)]_t+{\rm div}[{\rm div}(\r u\otimes u)]=\Delta F.
\end{equation}
Consider the following three elliptic problems on the whole space
$\mathbb{R}^2$:
\begin{equation}\label{xi}
-\Delta \x_1={\rm div}(\sqrt\r u(\sqrt\r-\sqrt{\bar\r})),
\end{equation}
\begin{equation}\label{xi2}
-\Delta \x_2=\sqrt{\bar\r} ~{\rm div}(\sqrt\r u),
\end{equation}
\begin{equation}\label{eta}
-\Delta \eta={\rm div}[{\rm div}(\r u\otimes u)],
\end{equation}
all with the boundary conditions $\x_1,\x_2,\eta\rightarrow 0$ as
$|x|\rightarrow\i$.

By the elliptic estimates and H${\rm\ddot{o}}$lder inequality, it
holds that
\begin{Lemma}\label{lemma4}
\begin{itemize}
\item[(1)] $\|\nabla\x_1\|_{2m}\leq Cm\|\r-\bar\r\|_{\f{2mk}{k-1}}\|u\|_{2mk},$
for any $k>1,m\geq1;$
\item[(2)] $\|\nabla\x_2\|_{2m}\leq Cm\Big[\|\r-\bar\r\|_{\f{2mk}{k-1}}\|u\|_{2mk}+\sqrt{\bar\r}\|u\|_{2m}\Big],$
for any $k>1,m\geq1;$
\item[(3)] $\|\nabla\x_2|x|^{\f\a2}\|_{2}\leq C\|\sqrt\r u|x|^{\f\a2}\|_{2},$
for $\a$ satisfying \eqref{alpha};
\item[(4)] $\|\eta\|_{2m}\leq Cm\big[\|\r-\bar\r\|_{\f{2mk}{k-1}}\|u\|^2_{4mk}+\bar\r\|u\|_{4m}^2\big],$
for any $k>1,m\geq1;$
\end{itemize} where $C$ are positive constants independent of $m,k$
and $r$.
\end{Lemma}
{\bf Proof:} By the elliptic estimates to the equations \eqref{xi},
\eqref{xi2}, respectively, and then using the H${\rm\ddot{o}}$lder
inequality, one has for any $k>1,m\geq1$,
$$
\begin{array}{ll}
\|\nabla\x_1\|_{2m}\leq C m\|\sqrt\r
u(\sqrt\r-\sqrt{\bar\r})\|_{2m}=
C m\| u(\r-\bar\r)\f{\sqrt\r}{\sqrt\r+\sqrt{\bar\r}}\|_{2m}\\
\di \qquad \leq
Cm\|\r-\bar\r\|_{\f{2mk}{k-1}}\|\f{\sqrt\r}{\sqrt\r+\sqrt{\bar\r}}\|_\i\|u\|_{2mk}\leq
Cm\|\r-\bar\r\|_{\f{2mk}{k-1}}\|u\|_{2mk},
\end{array}
$$
and
$$
\begin{array}{ll}
\|\nabla\x_2\|_{2m}\leq C m\|\sqrt\r u\|_{2m}\leq
C m\big[\| (\sqrt\r-\sqrt{\bar\r})u\|_{2m}+\sqrt{\bar\r}\|u\|_{2m}\big]\\
\di \qquad \leq
Cm\big[\|\sqrt\r-\sqrt{\bar\r}\|_{\f{2mk}{k-1}}\|u\|_{2mk}+\sqrt{\bar\r}\|u\|_{2m}\big]\leq
Cm\big[\|\r-\bar\r\|_{\f{2mk}{k-1}}\|u\|_{2mk}+\sqrt{\bar\r}\|u\|_{2m}\big].
\end{array}
$$
Thus the proofs of (1) and (2) are completed.

 By the similar proof as in Lemma 2.3 (2) in \cite{JWX3}, the
statements (3) can be proved.

Now we prove (4). By the elliptic estimates to the equation
\eqref{eta} and then using the H${\rm\ddot{o}}$lder inequality, one
has for any $k>1,m\geq1$,
$$
\begin{array}{ll}
\|\eta\|_{2m}\leq C m\|\r |u|^2\|_{2m}=
C m\big[\| (\r-\bar\r)|u|^2\|_{2m}+\bar\r\||u|^2\|_{2m}\big]\\
\di \qquad \leq
Cm\big[\|\r-\bar\r\|_{\f{2mk}{k-1}}\|u\|^2_{4mk}+\bar\r\|u\|_{4m}^2\big].
\end{array}
$$
Thus Lemma  \ref{lemma4} is proved.
 $\hfill\Box$

Based on Lemmas \ref{lemma1}-\ref{lemma3} and Lemma \ref{lemma4}, it
holds that

\begin{Lemma}\label{lemma5}
\begin{itemize}
\item[(1)] $\|\x_1\|_{2m}\leq Cm^{\f12}\|\r-\bar\r\|_{2m},$
for any $m\geq2;$
\item[(2)] $\|\x_2\|_{2m}\leq Cm^{\f12}\|\sqrt\r u|x|^{\f\a2}\|_{2}^{\f{2}{m\a}},$
for any $m+1\geq\f{4}{\a}$ and $\a$ satisfying \eqref{alpha};
\item[(3)] $\|u\|_{2m}\leq C m^{\f12}\big[\|\nabla u\|_{2}+1\big],$
for any $m\geq1;$
\item[(4)] $\|\nabla\x_1\|_{2m}\leq C m^{\f32}k^{\f12}\|\r-\bar\r\|_{\f{2mk}{k-1}}(\|\nabla u\|_2+1),$
for any $k>1, m\geq1;$
\item[(5)] $\|\nabla\x_2\|_{2m}\leq C m^{\f32}\big[k^{\f12}\|\r-\bar\r\|_{\f{2mk}{k-1}}+1\big](\|\nabla u\|_2+1),$
for any $k>1, m\geq1;$
\item[(6)] $\|\eta\|_{2m}\leq C m^2\big[k\|\r-\bar\r\|_{\f{2mk}{k-1}}+1\big]\big(\|\nabla u\|_2^2+1\big),$
for any $k>1, m\geq1;$
\end{itemize} where $C$ are positive constants independent of $m,k$.
\end{Lemma}
{\bf Proof:} (1) By Lemma \ref{lemma2}, it holds that
$$
\begin{array}{ll}
\di \|\x_1\|_{2m}\leq Cm^{\f12}\|\nabla\x_1\|_{\f{2m}{m+1}}&\di\leq
Cm^{\f12}\|\sqrt\r u\|_2\|\sqrt\r-\sqrt{\bar\r}\|_{2m}\leq Cm^{\f12}\|\r-\bar\r\|_{2m},
\end{array}
$$
where in the last inequality one has used the elementary energy
estimates.

(2). If $m+1>\f{4}{\a}$, then by interpolation inequality,
Caffarelli-Kohn-Nirenberg inequality and Lemma \ref{lemma4} (2), it
holds that
\begin{equation}\label{ine}
\|\x_2\|_{m+1}\leq \|\x_2\|^\t_{2m}\|\x_2\|_{\f{4}{\a}}^{1-\t}\leq
C\|\x_2\|^\t_{2m}\||x|^{\f\a2}\nabla \x_2\|_{2}^{1-\t}\leq
C\|\x_2\|^\t_{2m}\||x|^{\f\a2}\sqrt\r u\|_{2}^{1-\t}
\end{equation}
where
$$
\t=\f{\f{1}{m+1}-\f\a4}{\f{1}{2m}-\f\a4}.
$$
Then it follows from Lemma \ref{lemma1} (2) and \eqref{ine} that
$$
\begin{array}{ll}
\di \|\x_2\|_{2m}\leq Cm^{\f14}\|\nabla
\x_2\|_2^{\f12-\f{1}{2m}}\|\x_2\|_{m+1}^{\f12+\f{1}{2m}}\leq
Cm^{\f14}\|\sqrt\r
u\|_2^{\f12-\f{1}{2m}}\|\x_2\|_{2m}^{(\f12+\f{1}{2m})\t}\||x|^{\f\a2}\sqrt\r
u\|_{2}^{(\f12+\f{1}{2m})(1-\t)}\\
\di\qquad\quad  \leq
Cm^{\f14}\|\x_2\|_{2m}^{(\f12+\f{1}{2m})\t}\||x|^{\f\a2}\sqrt\r
u\|_{2}^{(\f12+\f{1}{2m})(1-\t)},
\end{array}
$$
which implies Lemma \ref{lemma5} (2) immediately.

Now we prove (3). First,
 $$
 \begin{array}{ll}
 \di \bar\r\int |u|^2dx=\int (\bar\r-\r)|u|^2dx+\int \r |u|^2dx\\
 \di =\int (\bar\r-\r)\big({\bf
1}|_{\{0\leq\r\leq2\bar\r\}}+{\bf 1}|_{\{\r>2\bar\r\}}\big)|u|^2dx+\int \r |u|^2dx\\
\di\leq \|(\bar\r-\r){\bf
1}|_{\{0\leq\r\leq2\bar\r\}}\|_2\|u\|_4^2+\|(\bar\r-\r){\bf 1}|_{\{\r>2\bar\r\}}\|_\g\|u\|_{\f{2\g}{\g-1}}^2+C\\
\di \leq \|\Psi(\r,\bar\r)\|_1^{\f12}\|u\|_2\|\nabla u\|_2+\|\Psi(\r,\bar\r)\|_1^{\f1\g}\|u\|_2^{2-\f2\g}\|\nabla u\|_2^{\f2\g}+C\\
\di \leq \s\|u\|_2^2+C_\s\|\nabla u\|_2^2+C.
 \end{array}
$$
Choosing $\s=\f{\bar\r}{2}$ in the above inequality yields that
\begin{equation}\label{u}
\|u\|_2^2\leq C\big(\|\nabla u\|_2^2+1\big).
\end{equation}
By Lemma \ref{lemma1} and the interpolation inequality, it holds that
$$
\|u\|_{2m}\leq C m^{\f14}\|\nabla u\|_2^{\f12-\f1{2m}}\|u\|_{m+1}^{\f12+\f1{2m}}\leq C m^{\f14}\|\nabla u\|_2^{\f12-\f1{2m}}\big(\|u\|^{\f{1}{m+1}}_2\|u\|^{\f m{m+1}}_{2m}\big)^{\f12+\f1{2m}},
$$
thus one has
$$
\|u\|_{2m}\leq C m^{\f12}\|\nabla u\|_2^{1-\f1m}\|u\|_2^{\f1m}\leq C m^{\f12}(\|\nabla u\|_2+1),
$$
where in the last inequality we have used \eqref{u}. The statement (3) is proved.

The assertions (3), (4) and (5) in Lemma \ref{lemma5} are the direct
consequences of Lemma \ref{lemma5} (2) and Lemma \ref{lemma4} (1),
(2), (4), respectively. Thus the proof of Lemma \ref{lemma5} is
completed. $\hfill\Box$

Substituting \eqref{xi}, \eqref{xi2} and \eqref{eta} into \eqref{vis-f} yields
that
$$
-\Delta\Big(\x_{1t}+\xi_{2t}+\eta+F\Big)=0,
$$
which implies that
\begin{equation*}
\x_{1t}+\x_{2t}+\eta+F=0.
\end{equation*}
It follows from the definition \eqref{flux} of the effective viscous
flux $F$ that
\begin{equation*}
\x_{1t}+\x_{2t}+(2\mu+\l(\r)){\rm div}u-(P(\r)-P(\bar\r))+\eta=0.
\end{equation*}
Then the continuity equation $\eqref{CNS}_1$ yields that
\begin{equation*}
\x_{1t}+\x_{2t}-\f{2\mu+\l(\r)}{\r}(\r_t+u\cdot\nabla \r)-(P(\r)-P(\bar\r))+\eta=0.
\end{equation*}
Define
\begin{equation*}
\Lambda(\r)=\int_{\bar\r}^\r\f{2\mu+\l(s)}{s}ds=2\mu\ln\f{\r}{\bar\r}+\f{1}{\b}(\r^\b-\bar\r^\b).
\end{equation*}
Then we obtain a new transport equation
\begin{equation}\label{transport-e}
(\Lambda(\r)-\x_1-\x_2)_t+u\cdot\nabla(\Lambda(\r)-\x_1-\x_2)+(P(\r)-P(\bar\r))+u\cdot\nabla(\x_1+\x_2)-\eta=0,
\end{equation}
which is crucial in the following Lemma for the higher integrability
of the density function.
\begin{Lemma}\label{lemma-rho}
  For any $k\geq2$ and $\b>1$, it holds that
  \begin{equation}\label{density-e}
    \sup_{t\in[0,T]}\|(\r-\bar\r)(t,\cdot)\|_k\leq C k^{\f{2}{\b-1}}.
  \end{equation}
\end{Lemma}
{\bf Proof:} Multiplying the equation \eqref{transport-e} by
$\r[(\Lambda(\r)-\x_1-\x_2)_+]^{2m-1}$ with $m$ being sufficiently large integer,
here and in what follows, the notation $(\cdots)_+$ denotes the
positive part of $(\cdots)$, one can get that
\begin{equation}\label{trans-e1}
\begin{array}{ll}
\di \f{1}{2m}\f{d}{dt}\int\r[(\Lambda(\r)-\x_1-\x_2)_+]^{2m}dx+\int\r
(P(\r)-P(\bar\r)){\bf
1}|_{\{\r>2\bar\r\}}[(\Lambda(\r)-\x_1-\x_2)_+]^{2m-1}dx\\
\di =-\int\r (P(\r)-P(\bar\r)){\bf
1}|_{\{0\leq\r\leq2\bar\r\}}[(\Lambda(\r)-\x_1-\x_2)_+]^{2m-1}dx\\
\di \quad+\int\r \eta[(\Lambda(\r)-\x_1-\x_2)_+]^{2m-1}dx -\int\r
u\cdot\nabla (\x_1+\x_2)[(\Lambda(\r)-\x_1-\x_2)_+]^{2m-1}dx.
\end{array}
\end{equation}
Denote
\begin{equation*}
f(t)=\big\{\int\r[(\Lambda(\r)-\x_1-\x_2)_+]^{2m} dx\big\}^{\f{1}{2m}},\qquad
t\in[0,T].
\end{equation*}
Now we estimate the three terms on the right hand side of
\eqref{trans-e1}. First, it holds that
\begin{equation}\label{rho-0}
\begin{array}{ll}
\di |\int\r (P(\r)-P(\bar\r)){\bf 1}|_{\{0\leq\r\leq2\bar\r\}}[(\Lambda(\r)-\x_1-\x_2)_+]^{2m-1}dx|\\
 \di \leq f(t)^{2m-1}\Big(\int \r|P(\r)-P(\bar\r)|^{2m}{\bf
 1}|_{\{0\leq\r\leq2\bar\r\}}dx\Big)^{\f1{2m}}\\
 \di \leq C f(t)^{2m-1}\|(\r-\bar\r){\bf
1}|_{\{0\leq\r\leq2\bar\r\}}\|_{2}\leq C
 f(t)^{2m-1}.
\end{array}
\end{equation}
Then, it follows that
\begin{equation}\label{rho-1}
\begin{array}{ll}
\di |\int\r \eta[(\Lambda(\r)-\x_1-\x_2)_+]^{2m-1}dx|\leq
\int\r^{\f{1}{2m}}|\eta|\big[\r(\Lambda(\r)-\x_1-\x_2)^{2m}_+\big]^{\f{2m-1}{2m}}dx\\
\di =\int\big[(\r-\bar\r)+\bar\r\big]^{\f{1}{2m}}|\eta|\big[\r(\Lambda(\r)-\x_1-\x_2)^{2m}_+\big]^{\f{2m-1}{2m}}dx\\
\di\leq
C\Big[\|\r-\bar\r\|_{2m\b+1}^{\f{1}{2m}}\|\eta\|_{2m+\f{1}{\b}}+\|\eta\|_{2m}\Big]\|\r(\Lambda(\r)-\x_1-\x_2)^{2m}_+\|_1^{\f{2m-1}{2m}}\\
\di\leq
C\Big[\|\r-\bar\r\|_{2m\b+1}^{\f{1}{2m}}(m+\f{1}{2\b})^2\big(k_1\|\r-\bar\r\|_{\f{2(m+\f{1}{2\b})k_1}{k_1-1}}+1\big)\\
\di\qquad\qquad\qquad\qquad  +m^2\big(k_2\|\r-\bar\r\|_{\f{2mk_2}{k_2-1}}+1\big)\Big]\big(\|\nabla u\|_2^2+1\big)f(t)^{2m-1}\\
\di\leq
Cm^2\big[\|\r-\bar\r\|_{2m\b+1}^{1+\f{1}{2m}}+1\big]\big(\|\nabla
u\|_2^2+1\big)f(t)^{2m-1},
\end{array}
\end{equation}
where in the last inequality we have taken $k_1=\f{\b}{\b-1}$ and
$k_2=\f{2m\b+1}{2m(\b-1)+1}.$

Next, for $\f{1}{2m\b+1}+\f1{p_1}+\f1{q_1}=1$ and
$\f1{p_2}+\f1{q_2}=1$ with $p_i,q_i> 1,(i=1,2)$, one has
\begin{equation}\label{rho-2}
\begin{array}{ll}
\di |-\int\r
u\cdot\nabla(\x_1+\x_2)[(\Lambda(\r)-\x_1-\x_2)_+]^{2m-1}dx|\\
\qquad\di \leq
\int\big[(\r-\bar\r)+\bar\r\big]^{\f{1}{2m}}|u||\nabla(\x_1+\x_2)|\big[\r(\Lambda(\r)-\x_1-\x_2)^{2m}_+\big]^{\f{2m-1}{2m}}dx\\
\qquad\di\leq
C\Big[\|\r-\bar\r\|_{2m\b+1}^{\f{1}{2m}}\|u\|_{2mp_1}\|\nabla(\x_1+\x_2)\|_{2mq_1}+\|u\|_{2mp_2}\|\nabla(\x_1+\x_2)\|_{2mq_2}\Big]f(t)^{2m-1}\\
\qquad\di\leq
C\Big[\|\r-\bar\r\|_{2m\b+1}^{\f{1}{2m}}(mp_1)^{\f12}(mq_1)^{\f32}\big(k_1^{\f12}\|\r-\bar\r\|_{\f{2mq_1k_1}{k_1-1}}+1\big)
  \\
  \qquad\qquad \di +(mp_2)^{\f12}(mq_2)^{\f32}\big(k_2\|\r-\bar\r\|_{\f{2mq_2k_2}{k_2-1}}+1\big)
\Big]\big(\|\nabla u\|_2^2+1\big)f(t)^{2m-1}\\
\qquad\di\leq
Cm^2\big[\|\r-\bar\r\|_{2m\b+1}^{1+\f{1}{2m}}+1\big]\big(\|\nabla
u\|_2^2+1\big)f(t)^{2m-1},
\end{array}
\end{equation}
where in the last inequality one has chosen
$p_1=\f{(2m\b+1)(\b+1)}{2m\b(\b-1)},q_1=\f{(\b+1)(2m\b+1)}{4m\b}, k_1=\f{2\b}{\b-1}$ and $p_2=\f{2\b}{\b-1}, q_2=\f{2\b}{\b+1},
k_2=\f{(\b+1)(2m\b+1)}{2m\b(\b-1)+(\b+1)}.$

Substituting \eqref{rho-0}, \eqref{rho-1} and \eqref{rho-2} into
\eqref{trans-e1} yields that
\begin{equation*}
\begin{array}{ll}
\di\f{1}{2m}\f{d}{dt}(f^{2m}(t))+\int\r (P(\r)-P(\bar\r)){\bf
1}|_{\{\r>2\bar\r\}}[(\Lambda(\r)-\x_1-\x_2)_+]^{2m-1}dx\\
\di\leq
Cm^2\big[\|\r-\bar\r\|_{2m\b+1}^{1+\f{1}{2m}}+1\big]\big(\|\nabla
u\|_2^2+1\big)f(t)^{2m-1}.
\end{array}
\end{equation*}
Then it holds that
\begin{equation*}
\f{d}{dt}f(t)\leq
Cm^2\big[\|\r-\bar\r\|_{2m\b+1}^{1+\f{1}{2m}}+1\big]\big(\|\nabla
u\|_2^2+1\big).
\end{equation*}
Integrating the above inequality over $[0,t]$ gives that
\begin{equation}\label{f}
f(t)\leq
f(0)+Cm^2\int_0^t\big[\|\r-\bar\r\|_{2m\b+1}^{1+\f{1}{2m}}+1\big]\big(\|\nabla
u\|_2^2+1\big)d\tau.
\end{equation}
Now we calculate the quantity
$$f(0)=\Big(\int\r_0[(\Lambda(\r_0)-\x_{10}-\x_{20})_+]^{2m} dx\Big)^{\f{1}{2m}}.$$
By Lemma \ref{lemma4} (1), (2) and Lemma \ref{lemma5} (1), (2) with
$t=0$, we can easily get
$$
\|\x_{10}+\x_{20}\|_{L^\i}\leq C.
$$
Furthermore, by the definition of
$\Lambda(\r_0)=2\mu\ln\f{\r_0}{\bar\r}+\f1\b((\r_0)^\b-\bar\r^\b)$, we have
$$
\Lambda(\r_0)-\x_{10}-\x_{20}\rightarrow-\i, \quad{\rm as}\quad \r_0\rightarrow0+.
$$
Thus there exists a positive constant $\s_0$, such that if
$0\leq\r_0\leq\s_0$, then
$$
(\Lambda(\r_0)-\x_{10}-\x_{20})_+\equiv0.
$$
Now one has
\begin{equation}\label{f0}
\begin{array}{ll}
f(0)&\di =\Big[\Big(\int_{[0\leq\r_0\leq
\s_0]}+\int_{[\s_0\leq\r_0\leq
M]}\Big)\r_0(\Lambda(\r_0)-\x_{10}-\x_{20})_+^{2m}
dx\Big]^{\f{1}{2m}}\\
&\di=\Big[\int_{[\s_0\leq\r_0\leq
M]}\r_0(\Lambda(\r_0)-\x_{10}-\x_{20})_+^{2m} dx\Big]^{\f{1}{2m}}\\
&\di \leq C(\s_0,M)\Big[\|(\r_0-\bar\r){{\bf 1}_{\s_0\leq\r_0\leq
M}}\|_{2m}+\|\x_{10}+\x_{20}\|_{2m}\Big]\leq C(\s_0,M) m^{\f32},
\end{array}
\end{equation}
where the positive constant $C(\s_0,M)$ is independent of $m$ and the lower bound of the density.

It follows from \eqref{f} and \eqref{f0} that for $t\in[0,T]$,
\begin{equation}\label{f1}
f(t)\leq
C m^2\Big[1+\int_0^t\big(\|\r-\bar\r\|_{2m\b+1}^{1+\f{1}{2m}}+1\big)\big(\|\nabla
u\|_2^2+1\big)d\tau\Big].
\end{equation}
For any $t\in[0,T]$,  set
$\Omega_1(t)=\{x\in\mathbb{R}^2|\r(t,x)>2\bar\r\}$ and
$\Omega_2(t)=\{x\in\Omega_1(t)|(\Lambda(\r)-\xi_1-\x_2)(t,x)>0\}$. Then one has on $\Omega_1(t),$ $|\r-\bar\r|^\b\leq C\beta|\Lambda(\r)|$ for some  constant $C>0$, and on $\Omega_1(t)\setminus\Omega_2(t),$ $0<\Lambda(\r)\leq \x_1+\x_2.$  Thus it
holds that
\begin{equation}\label{rs}
\begin{array}{ll}
\di\|\r-\bar\r\|_{2m\b+1}^\b(t)=\Big(\int_{\Omega_1(t)}|\r-\bar\r|^{2m\b+1}dx+\int_{\mathbb{R}^2\setminus\Omega_1(t)}|\r-\bar\r|^{2m\b+1}dx\Big)^{\f{\b}{2m\b+1}}\\
\di\leq\Big(\int_{\Omega_1(t)}|\r-\bar\r|^{2m\b+1}dx+\bar\r^{2m\b-1}\int_{\mathbb{R}^2\setminus\Omega_1(t)} |\r-\bar\r|^2 dx\Big)^{\f{\b}{2m\b+1}}\\
\di \leq
\Big[\int_{\Omega_1(t)}\big(|\r-\bar\r|^\b\big)^{\f{2m\b+1}{\b}}dx\Big]^{\f\b{2m\b+1}}+C
\leq
C\Big(\int_{\Omega_1(t)}|\b\Lambda(\r)|^{\f{2m\b+1}{\b}}dx\Big)^{\f{\b}{2m\b+1}}+C\\
\di \leq
C\Big(\int_{\Omega_1(t)}\Lambda(\r)^{2m}\Lambda(\r)^{\f{1}{\b}}dx\Big)^{\f{\b}{2m\b+1}}+C\leq
C\Big(\int_{\Omega_1(t)}\r\Lambda(\r)^{2m}dx\Big)^{\f{\b}{2m\b+1}}+C\\
\di=
C\Big(\int_{\Omega_2(t)}\r|\Lambda(\r)-\x_1-\x_2+(\x_1+\x_2)|^{2m}dx+\int_{\Omega_1(t)\setminus\Omega_2(t)}\r|\Lambda(\r)|^{2m}dx\Big)^{\f{\b}{2m\b+1}}+C\\
\di\leq
C\Big[\int_{_{\Omega_2(t)}}\r(\Lambda(\r)-\x_1-\x_2)^{2m}dx+\int_{_{\Omega_2(t)}}\r|\x_1+\x_2|^{2m}dx+\int_{_{\Omega_1(t)\setminus\Omega_2(t)}}\r|\x_1+\x_2|^{2m}dx\Big]^{\f{\b}{2m\b+1}}+C\\
\di\leq
C\Big(f(t)^{2m}+\int_{\mathbb{R}^2}\r|\x_1+\x_2|^{2m}dx\Big)^{\f{\b}{2m\b+1}}+C\\
\di\leq
C\Big[f(t)+\Big(\int_{\mathbb{R}^2}\r|\x_1|^{2m}dx\Big)^{\f{\b}{2m\b+1}}+\Big(\int_{\mathbb{R}^2}\r|\x_2|^{2m}dx\Big)^{\f{\b}{2m\b+1}}+1\Big].
\end{array}
\end{equation}
Note that
\begin{equation}\label{rs1}
\begin{array}{ll}
\di\Big(\int_{\mathbb{R}^2}\r|\x_1|^{2m}dx\Big)^{\f{\b}{2m\b+1}}&\di
\leq
C\Big(\int_{\mathbb{R}^2}|\r-\bar\r||\x_1|^{2m}dx\Big)^{\f{\b}{2m\b+1}}+C\Big(\int_{\mathbb{R}^2}|\x_1|^{2m}dx\Big)^{\f{\b}{2m\b+1}}:=K_{11}+K_{12}.
\end{array}
\end{equation}
\begin{equation}\label{K11}
\begin{array}{ll}
K_{11}&\di \leq C\|\r-\bar\r\|_{2m\b+1}^{\f{\b}{2m\b+1}}\||\x_1|^{2m}\|^{\f{\b}{2m\b+1}}_{\f{2m\b+1}{2m\b}}=C\|\r-\bar\r\|_{2m\b+1}^{\f{\b}{2m\b+1}}\|\x_1\|^{\f{2m\b}{2m\b+1}}_{2m+\f{1}{\b}}\\
&\di \leq\|\r-\bar\r\|_{2m\b+1}^{\f{\b}{2m\b+1}}\Big[C(m+\f{1}{2\b})^{\f12}\|\r-\bar\r\|_{2m+\f{1}{\b}}\Big]^{\f{2m\b}{2m\b+1}}\\
&\di \leq
Cm^{\f12}\|\r-\bar\r\|_{2m\b+1}^{\f{\b}{2m\b+1}}\Big[\|(\r-\bar\r){\bf 1}|_{\{0\leq\r\leq2\bar\r\}}\|^{\f{2m\b}{2m\b+1}}_2\\
&\di\qquad\qquad\qquad\qquad\qquad  +\|(\r-\bar\r){\bf 1}|_{\{\r>2\bar\r\}}\|^{\f{2\g m\b(\b-1)}{(2m\b-\g+1)(2m\b+1)}}_\g\|\r-\bar\r\|^{\f{2m\b(2m\b-\g\b+1)}{(2m\b-\g+1)(2m\b+1)}}_{2m\b+1}\Big]\\
&\di \leq Cm^{\f12}\big[\|\r-\bar\r\|_{2m\b+1}+1\big],
\end{array}
\end{equation}
and
\begin{equation}\label{K12}
K_{12} =\|\x_1\|^{\f{2m\b}{2m\b+1}}_{2m} \leq\Big(Cm^{\f12}\|\r-\bar\r\|_{2m}\Big)^{\f{2m\b}{2m\b+1}}
\leq Cm^{\f12}\big[\|\r-\bar\r\|_{2m\b+1}+1\big].
\end{equation}
Furthermore, it holds that
\begin{equation}\label{rs11}
\begin{array}{ll}
\di\Big(\int_{\mathbb{R}^2}\r|\x_2|^{2m}dx\Big)^{\f{\b}{2m\b+1}}
\leq C\Big(\int_{\mathbb{R}^2}|\r-\bar\r||\x_2|^{2m}dx\Big)^{\f{\b}{2m\b+1}}+C\Big(\int_{\mathbb{R}^2}|\x_2|^{2m}dx\Big)^{\f{\b}{2m\b+1}}\\
\di\qquad\qquad \leq C\|\r-\bar\r\|_{2m\b+1}^{\f{\b}{2m\b+1}}\|\x_2\|^{\f{2m\b}{2m\b+1}}_{2m+\f1\b}+C\|\x_2\|^{\f{2m\b}{2m\b+1}}_{2m}\\
\di\qquad\qquad\leq C m^{\f12}\|\r-\bar\r\|_{2m\b+1}^{\f{\b}{2m\b+1}}\|\sqrt\r u|x|^{\f\a2}\|_2^{\f{2}{(m+\f{1}{2\b})\a}\f{2m\b}{2m\b+1}}+Cm^{\f12}\|\sqrt\r u|x|^{\f\a2}\|_2^{\f{2}{m\a}\f{2m\b}{2m\b+1}}\\
\di \qquad\qquad \leq C
\big[\|\sqrt\r
u|x|^{\f\a2}\|_2^2+m^{\f12}\|\r-\bar\r\|_{2m\b+1}+m^\f12\big].
\end{array}
\end{equation}
Substituting and \eqref{rs1}, \eqref{K11}, \eqref{K12} and \eqref{rs11}  into \eqref{rs} yields that
\begin{equation}\label{rs2}
\begin{array}{ll}
\|\r-\bar\r\|_{2m\b+1}^\b(t)&\di\leq
C\Big[m^{\f12}+f(t)+m^{\f12}\|\r-\bar\r\|_{2m\b+1}(t)+\|\sqrt\r u|x|^{\f\a2}\|_2^2(t)\Big]\\
&\di \leq
\f12\|\r-\bar\r\|_{2m\b+1}^\b(t)+C\Big[f(t)+m^{\f{\b}{2(\b-1)}}+\|\sqrt\r
u|x|^{\f\a2}\|^2_2(t)\Big].
\end{array}
\end{equation}
Thus it follows from \eqref{f1}, \eqref{rs2} and the weighted estimates in Lemma
\ref{lemma-wee} that
\begin{equation*}
\begin{array}{ll}
\|\r-\bar\r\|_{2m\b+1}^\b(t)\di\leq C\Big[f(t)+m^{\f{\b}{2(\b-1)}}+\|\sqrt\r u|x|^{\f\a2}\|_2^2(t)\Big]\\
\di\leq
C\Big[m^2+m^2\int_0^t\|\r-\bar\r\|_{2m\b+1}^{1+\f{1}{2m}}\big(\|\nabla
u\|_2^2+1\big)d\tau
+\int_0^t\|\r-\bar\r\|_{2m\b+1}^{\b}\big(\|\nabla
u\|_2^2+1\big)d\tau\Big].
\end{array}
\end{equation*}
Applying Gronwall's inequality to the above inequality yields that
\begin{equation*}
\|\r-\bar\r\|_{2m\b+1}^\b(t)\leq
C m^2 \Big[1+\int_0^t\|\r-\bar\r\|_{2m\b+1}^{1+\f{1}{2m}}\big(\|\nabla
u\|_2^2(\tau)+1\big)d\tau\Big].
\end{equation*}
Denote
$$
y(t)=m^{-\f{2}{\b-1}}\|\r-\bar\r\|_{2m\b+1}(t).
$$
Then it holds that
$$
\di y^\b(t)\leq
C\Big[1+\int_0^ty(\tau)^{1+\f{1}{2m}}\|\nabla u\|_2^2(\tau)d\tau\Big] \leq C\Big[1+\int_0^t\big(y^\b(\tau)+1\big)\|\nabla
u\|_2^2(\tau) d\tau\Big].
$$
So applying the Gronwall's inequality to the above inequality yields
that
$$
y(t)\leq C,\quad \forall t\in[0,T],
$$
that is, for sufficiently large $m>1$,
$$
\|\r-\bar\r\|_{2m\b+1}(t)\leq Cm^{\f{2}{\b-1}},\quad \forall t\in[0,T].
$$
Equivalently,  \eqref{density-e} holds for sufficiently large $k$. Now by the elementary energy estimate Lemma \ref{lemma-ee}, if $\g\geq2,$ then
\begin{equation}\label{g1}
\|\r-\bar\r\|_2(t)\leq C\|\Psi(\r,\bar\r)\|_1^\f12\leq C,
\end{equation}
and if $1<\g<2$, then
\begin{equation}\label{g2}
\|(\r-\bar\r){\bf 1}|_{\{0\leq\r\leq2\bar\r\}}\|_2(t)\leq C\|\Psi(\r,\bar\r){\bf 1}|_{\{0\leq\r\leq2\bar\r\}}\|_1^\f12\leq C,
\end{equation}
and
$$
\|(\r-\bar\r){\bf 1}|_{\{\r>2\bar\r\}}\|_\g(t)\leq C\|\Psi(\r,\bar\r){\bf 1}|_{\{\r>2\bar\r\}}\|_1^\f1\g\leq C.
$$
Therefore, for $1<\g<2$, it holds that
\begin{equation}\label{g3}
\|(\r-\bar\r){\bf 1}|_{\{\r>2\bar\r\}}\|_2(t)\leq\|(\r-\bar\r){\bf 1}|_{\{\r>2\bar\r\}}\|_\g^\t\|\r-\bar\r\|^{1-\t}_{k}\leq C,
\end{equation}
where $k$ is sufficiently large such that \eqref{density-e} holds and $\t\in(0,1)$ satisfying $\f12=\f{\t}{\g}+\f{1-\t}{k}.$ Thus by \eqref{g1}, \eqref{g2} and \eqref{g3}, it holds that for any $\g>1$ and $t\in[0,T]$,
\begin{equation*}
\|\r-\bar\r\|_2(t)\leq C.
\end{equation*}
Thus Lemma \ref{lemma-rho}
is proved for any $k\geq2$. $\hfill\Box$

\underline{Step 4: First-order derivative estimates of the
velocity.}

Set
$$
Z^2(t)=\int(\mu\omega^2+\f{F^2}{2\mu+\l(\r)}) dx,
$$
$$
\varphi^2(t)=\int\r(H^2+L^2)dx=\int\r |\dot u|^2 dx
$$
and
$$
\Phi_T=\sup_{t\in[0,T]}\|\r(\cdot,t)\|_\i+1.
$$
The following Lemma is motivated by \cite{P}.

\begin{Lemma}\label{lemma-u-der}
For any $\v>0$, there exists a positive constant $C_\v$, such that
  \begin{equation*}
    \sup_{t\in[0,T]}\log(e+Z^2(t))+\int_0^T\f{\varphi^2(t)}{e+Z^2(t)}dt\leq C_\v\Phi_T^{1+\v\b}.
  \end{equation*}
\end{Lemma}
{\bf Proof:} Multiplying the equation $\eqref{F-omega}_1$ by
$\mu\omega$, the equation $\eqref{F-omega}_2$ by
$\f{F}{2\mu+\l(\r)}$, respectively, and then summing the resulted
equations together, one has
\begin{equation}\label{ue1}
\begin{array}{ll}
\di\f12\f{d}{dt}\int(\mu\omega^2+\f{F^2}{2\mu+\l(\r)}) dx+\f{\mu}{2}\int\omega^2{\rm div }u dx-\f12\int\r F^2(\f{1}{2\mu+\l(\r)})^\prime{\rm div }u dx \\
\di  -\f12\int F^2\f{{\rm div}u}{2\mu+\l(\r)} dx-\int\r F({\rm div}
u)(\f{P(\r)-P(\bar\r)}{2\mu+\l(\r)})^\prime dx+\int
F[(u_{1x_1})^2+2u_{1x_2}u_{2x_1}+(u_{2x_2})^2]dx\\
\di =-\int\r(H^2+L^2) dx.
\end{array}
\end{equation}
Notice that
$$
\begin{array}{ll}
\di (u_{1x_1})^2+2u_{1x_2}u_{2x_1}+(u_{2x_2})^2\di =(u_{1x_1}+u_{2x_2})^2+2(u_{1x_2}u_{2x_1}-u_{1x_1}u_{2x_2})\\
\di\qquad =({\rm div} u)^2+2(u_{1x_2}u_{2x_1}-u_{1x_1}u_{2x_2})\\
\di\qquad =({\rm div}
u)\left(\f{F}{2\mu+\l(\r)}+\f{P(\r)-P(\bar\r)}{2\mu+\l(\r)}\right)+2(u_{1x_2}u_{2x_1}-u_{1x_1}u_{2x_2}),
\end{array}
$$
then one has
\begin{equation}\label{ue1}
\begin{array}{ll}
\di\f12\f{d}{dt}\int(\mu\omega^2+\f{F^2}{2\mu+\l(\r)}) dx+\int\r(H^2+L^2) dx\\
\di \quad =-\f{\mu}{2}\int\omega^2{\rm div }u dx+\f12\int F^2({\rm
div}u)\Big[\r(\f{1}{2\mu+\l(\r)})^\prime-\f{1}{2\mu+\l(\r)}\Big]
dx\\
\di \quad +\int F({\rm div}
u)\Big[\r(\f{P(\r)-P(\bar\r)}{2\mu+\l(\r)})^\prime-\f{P(\r)-P(\bar\r)}{2\mu+\l(\r)}\Big]
dx-\int 2F(u_{1x_2}u_{2x_1}-u_{1x_1}u_{2x_2})dx.
\end{array}
\end{equation}
Then
\begin{equation}\label{fact2}
\begin{array}{ll}
\di \|\nabla u\|_2+\|\omega\|_2+\|{\rm div}
u\|_2+\|(2\mu+\l(\r))^{\f12}{\rm div} u\|_2\\
\di \leq C\Big[Z(t)+\Big(\int\f{|P(\r)-P(\bar\r)|^2}{2\mu+\l(\r)}
dx\Big)^{\f12}\Big]\leq C(Z(t)+1).
\end{array}
\end{equation}

Now we estimate the four terms on the right hand side of
\eqref{ue1}. First, by the interpolation inequality, Lemma
\ref{lemma1} and \eqref{fact2}, it holds that
\begin{equation}\label{ue2}
\begin{array}{ll}
\di |-\f{\mu}{2}\int\omega^2{\rm div }u dx|\leq C\|{\rm div}
u\|_2\|\o\|_4^2 \leq
C(Z(t)+1)\|\o\|_2\|\nabla\o\|_{2}\\
\di \qquad\qquad\leq
C(Z(t)+1)\|\o\|_2\|\r\dot u\|_{2}\leq C(Z(t)+1)\|\r\|^{\f12}_\i\|\o\|_2\|\sqrt\r\dot u\|_{2}\\
\qquad\qquad\leq \s \|\sqrt\r\dot u\|_{2}^2+C_\s
(Z(t)^2+1)\|\r\|_\i\|\o\|_2^2.
\end{array}
\end{equation}
Next, one has
\begin{equation*}
\begin{array}{ll}
\di |\f12\int F^2{\rm
div}u\Big[\r(\f{1}{2\mu+\l(\r)})^\prime-\f{1}{2\mu+\l(\r)}\Big]
dx|\\
\di \leq C\int|F|^2
\f{|{\rm div} u|}{2\mu+\l(\r)}dx \leq
\|{\rm div} u\|_2\|\f{F^2}{2\mu+\l(\r)}\|_2,
\end{array}
\end{equation*}
while for any $\v>0$ suitably small,
\begin{equation}\label{ue3}
\begin{array}{ll}
\di
\|\f{F^2}{2\mu+\l(\r)}\|_2\leq \|\f{F}{\sqrt{2\mu+\l(\r)}}\|_2^{1-\v} \|F\|^{1+\v}_{\f{2(1+\v)}{\v}}\\
\di \leq C\|\f{F}{\sqrt{2\mu+\l(\r)}}\|_2^{1-\v} \big(\|F\|_2^{\f{\v}{1+\v}}\|\nabla F\|_2^{\f1{1+\v}}\big)^{1+\v}\\
\di \leq C\|\f{F}{\sqrt{2\mu+\l(\r)}}\|_2^{1-\v} \|F\|_2^{\v}\|\nabla F\|_2
\di \leq C\|\f{F}{\sqrt{2\mu+\l(\r)}}\|_2\|\r\|_\i^{\f{1+\b \v}{2}}\|\sqrt{\rho}\dot u\|_2.
\end{array}
\end{equation}
Then it holds that
\begin{equation}\label{ue40}
\begin{array}{ll}
\di |\f12\int F^2{\rm
div}u\Big[\r(\f{1}{2\mu+\l(\r)})^\prime-\f{1}{2\mu+\l(\r)}\Big]
dx|\\
\di \leq C
\|{\rm div} u\|_2\|\f{F}{\sqrt{2\mu+\l(\r)}}\|_2\|\r\|_\i^{\f{1+\b \v}{2}}\|\sqrt{\rho}\dot u\|_2\\
\di \leq \sigma \|\sqrt{\rho}\dot u\|_2^2+C_\s \|\r\|_\i^{1+\b \v}\|{\rm div} u\|^2_2\|\f{F}{\sqrt{2\mu+\l(\r)}}\|^2_2.
\end{array}
\end{equation}
On the other hand, it holds that
\begin{equation}\label{ue4}
\begin{array}{ll}
\di |\int F({\rm div}
u)\Big[\r(\f{P(\r)-P(\bar\r)}{2\mu+\l(\r)})^\prime-\f{P(\r)-P(\bar\r)}{2\mu+\l(\r)}\Big]
dx|\\
\di \leq C \int |F||{\rm div}
u|\f{|P(\r)-P(\bar\r)|+1}{2\mu+\l(\r)} dx\\
\di \leq C \|{\rm div} u\|_2\Big[\|F\|_{\f{2(2+\v)}{\v}}\|P(\r)-P(\bar\r)\|_{2+\v}+\|\f{F}{\sqrt{2\mu+\l(\r)}}\|_2\Big]\\
\di \leq C \|{\rm div} u\|_2\Big[\|\f{F}{\sqrt{2\mu+\l(\r)}}\|_2^{\f{\v}{2+\v}}\|\r\|_\i^{\f{\b\v}{2(2+\v)}}\|\nabla F\|_2^{\f{2}{2+\v}}+\|\f{F}{\sqrt{2\mu+\l(\r)}}\|_2\Big]\\
\di \leq C \|{\rm div} u\|_2\Big[\|\f{F}{\sqrt{2\mu+\l(\r)}}\|_2^{\f{\v}{2+\v}}\|\r\|_\i^{\f12+\f{\b\v}{2(2+\v)}}\|\sqrt\r \dot u\|_2^{\f{2}{2+\v}}+\|\f{F}{\sqrt{2\mu+\l(\r)}}\|_2\Big]\\
\di \leq \s \|\sqrt\r \dot u\|_2^2+C_\s\|\r\|_\i^{1+\f{\b\v}{2+\v}}\|{\rm div} u\|_2^{\f{2+\v}{1+\v}}\|\f{F}{\sqrt{2\mu+\l(\r)}}\|_2^{\f{\v}{1+\v}}+C \|{\rm div} u\|_2\|\f{F}{\sqrt{2\mu+\l(\r)}}\|_2\\
\di \leq \s \|\sqrt\r \dot u\|_2^2+C_\s(1+\|\r\|_\i)^{1+\b\v}\big(\|{\rm div} u\|_2^2+1\big)\big(\|\f{F}{\sqrt{2\mu+\l(\r)}}\|_2^2+1\big).
\end{array}
\end{equation}
Then due to \cite{P}, it holds that
\begin{equation}\label{ue5}
\begin{array}{ll}
\di |-\int 2F(u_{1x_2}u_{2x_1}-u_{1x_1}u_{2x_2})dx|=|-\int 2F \nabla
u_1\cdot \nabla^\perp u_2 dx|\\
\di \leq C\|F\|_{\rm BMO}\|\nabla
u_1\cdot \nabla^\perp u_2\|_{\mathcal{H}^1}\leq C\|\nabla F\|_{2}\|\nabla
u\|_2^2\leq C\|\r\|_\i^{\f12}\|\sqrt{\r} \dot u\|_2\|\nabla u\|_2^2\\
\di \leq \s \|\sqrt{\r} \dot u\|_2^2+C_\s \|\r\|_\i\|\nabla u\|_2^4\leq \s \|\sqrt{\r} \dot u\|_2^2+C_\s \|\r\|_\i\|\nabla u\|_2^2\Big[1+Z^2(t)\Big].
\end{array}
\end{equation}
In summary, substituting \eqref{ue2}, \eqref{ue40}, \eqref{ue4} and \eqref{ue5} into \eqref{ue1}, one can arrive at
\begin{equation*}
\di \f{d}{dt}Z^2(t)+\varphi^2(t)\leq C\Phi_T^{1+\b\v}(1+\|\nabla u\|_2^2)(1+Z^2(t))
\end{equation*}
Multiplying the above inequality by $\f{1}{e+Z^2(t)}$ and then integrating over $[0,T]$ give  the proof of Lemma \ref{lemma-u-der} . $\hfill\Box$

\underline{Step 5: Upper and lower bound of the density:}

The following Lemma comes from \cite{hl1, hl2}. With the following Lemma, the index $\b$ can be improved to $\b>\f43$ as in \cite{hl1, hl2}.

\begin{Lemma}\label{lemma-u-d}
  There exists a positive constant $C$, such that
  \begin{equation*}
\sup_{t\in[0,T]}\int\r |u|^{2+\nu} dx \leq C,
  \end{equation*}
  where $\nu=\f{\mu^{\f12}}{2(\mu+1)}\Phi_T^{-\f\b2}\in(0,\f14].$
\end{Lemma}

{\bf Proof:} Multiplying the momentum equation $\eqref{CNS}_2$ by $(2+\nu)u |u|^{\nu}$ and integrating over $\mathbb{R}^2$ with respect to $x$ lead to
$$
\begin{array}{ll}
\di \f{d}{dt}\int \r |u|^{2+\nu}dx+\mu(2+\nu) \int |\nabla u|^2|u|^{\nu}dx+(2+\nu)\int(\mu+\l(\r)) ({\rm div}u)^2 |u|^\nu dx\\
\di =(2+\nu)\int(P(\r)-P(\bar\r)){\rm div}(u|u|^\nu)dx-\mu(2+\nu)\int\nabla \f{|u|^2}{2}\cdot \nabla |u|^\nu dx\\
\di \quad -(2+\nu)\int (\mu+\l(\r))({\rm div} u) u\cdot \nabla|u|^\nu dx.
\end{array}
$$
Now we only estimate the first term on the right hand side of the above equality, since the other terms can be done similarly as in \cite{hl2}. Then it holds that
$$
\begin{array}{ll}
\di (2+\nu)|\int(P(\r)-P(\bar\r)){\rm div}(u|u|^\nu)dx|
\leq (2+\nu)(1+\nu) \int|P(\r)-P(\bar\r)||\nabla u||u|^\nu dx\\
\di \leq \s (2+\nu)\int |\nabla u|^2|u|^\nu dx+C_\s (2+\nu)(1+\nu)^2 \int|P(\r)-P(\bar\r)|^2 |u|^\nu dx\\
\di \leq \s (2+\nu)\int |\nabla u|^2|u|^\nu dx+C_\s (2+\nu)(1+\nu)^2 \|P(\r)-P(\bar\r)\|_{2q_1}^2 \|u\|_{q_2\nu}^\nu \\
\di \leq \s (2+\nu)\int |\nabla u|^2|u|^\nu dx+C_\s (2+\nu)(1+\nu)^2 \big(\|\nabla u\|_{2}^2+1\big).
\end{array}
$$
where $q_1, q_2>1$ satisfying $\f1{q_1}+\f1{q_2}=1.$ Thus Lemma \ref{lemma-u-d} is proved.$\hfill\Box$

Now one can obtain the upper and lower bound of the density by using the transport equation \eqref{transport-e}.

\begin{Lemma}\label{lemma-den}
  There exists positive constants $C_1$ and $c_1$ such that
  \begin{equation*}
   c_1\leq \r(t,x)\leq C_1, \qquad \forall (t,x)\in[0,T]\times\mathbb{R}^2.
  \end{equation*}
\end{Lemma}

{\bf Proof:}
First, for any $p>2$ and $q>1$ satisfying
$$
\f1p=\f{\f2p}{2+\nu}+\f{1-\f2p}{q},
$$
it holds that
$$
\begin{array}{ll}
\di
\|\r u\|_p\leq \|\r u\|_{2+\nu}^{\f{2}{p}}\|\r u\|_q^{1-\f2p}\leq C\Big(\|\r^\f{1}{2+\nu}u\|_{2+\nu}\|
\r\|_\i^{\f{1+\nu}{2+\nu}}\Big)^{\f2p} \big(\|\r\|_\i\|u\|_q\big)^{1-\f2p}\\
\di\qquad\quad  \leq C\|\r\|_\i^{1-\f2p+\f{2(1+\nu)}{p(2+\nu)}}\Big[q^\f12(\|\nabla u\|_2+1)\Big]^{1-\f2p},
\end{array}
$$
where in the last inequality one has used Lemma \ref{lemma5} (3). It can be computed that
$$
q=(1+\f 2\nu)(p-2)\leq C_p \Phi_T^{\f\b2}.
$$
Therefore, one has
\begin{equation}\label{ru1}
\|\r u\|_p\leq C\|\r\|_\i^{1-\f{2}{p(2+\nu)}}\Phi_T^{\f\b4(1-\f2p)}(\|\nabla u\|_2^{1-\f2p}+1)\leq C\Phi_T^{1+\f\b4}(\|\nabla u\|_2^{1-\f2p}+1).
\end{equation}
Note that by the definition of $\x_i~(i=1,2)$ from \eqref{xi} and \eqref{xi2}
\begin{equation}\label{ct}
u\cdot\nabla(\x_1+\x_2)-\eta=[u, R_iR_j](\r u),
\end{equation}
where $[\cdot,\cdot]$ is the usual commutator and $R_i, R_j$ are the Riesz operators. Thus from \eqref{transport-e}, it holds that
\begin{equation}\label{te1}
D_t\Lambda(\r)-D_t(\x_1+\x_2)+(P(\r)-P(\bar\r))+[u, R_iR_j](\r u)=0,
\end{equation}
where the material derivative $D_t:=\partial_t+u\cdot\nabla$.
Along the particle path $\vec{X}(\tau;t,x)$ through the point
$(t,x)\in[0,T]\times\mathbb{R}^2$ defined by
\begin{equation*}
\left\{
\begin{array}{ll}
\di \f{d\vec{X}(\tau;t,x)}{d\tau}=u(\tau,\vec{X}(\tau;t,x)),\\
 \di
\vec{X}(\tau;t,x)|_{\tau=t}=x,
\end{array}
\right.
\end{equation*}
from the equation $\eqref{te1}$, there holds the following ODE
\begin{equation*}
\begin{array}{ll}
\di \f{d}{d\tau} (\Lambda(\r)-\x_1-\x_2)(\tau,\vec{X}(\tau;t,x))+(P(\r)-P(\bar\r)){\bf
1}|_{\{\r>2\bar\r\}}(\tau,\vec{X}(\tau;t,x))\\
\di\qquad=-\Big((P(\r)-P(\bar\r)){\bf
1}|_{\{0\leq\r\leq2\bar\r\}}+[u, R_iR_j](\r u)\Big)(\tau,\vec{X}(\tau;t,x)),
\end{array}
\end{equation*}
and thus
\begin{equation*}
\f{d}{d\tau} (\Lambda(\r)-\x_1-\x_2)(\tau,\vec{X}(\tau;t,x))\leq -\Big((P(\r)-P(\bar\r)){\bf
1}|_{\{0\leq\r\leq2\bar\r\}}+[u, R_iR_j](\r u)\Big)(\tau,\vec{X}(\tau;t,x)).
\end{equation*}
Integrating the above inequality over $[0,t]$ yields that
\begin{equation}\label{theta-1}
\begin{array}{ll}
\di 2\mu\ln\f{\r(t,x)}{\r_0(\vec{X}_0)}+\f{1}{\b}\big(\r^\b(t,x)-\r_0^\b(\vec{X}_0)\big)-\big((\x_1+\x_2)(t,x)-(\x_{10}+\x_{20})(\vec{X}_0)\big)\\
\di \qquad\qquad\qquad\qquad\qquad \leq C+\int_0^t\|[u, R_iR_j](\r u)\|_\i ds,
\end{array}
\end{equation}
with $\vec{X}_0=\vec{X}(\tau;t,x)|_{\tau=0}$.

 Then for any sufficiently large $p>4$, by the commutator estimates for \eqref{ct} and \eqref{ru1}, it holds that
 $$
 \begin{array}{ll}
 \di \|[u, R_iR_j](\r u)\|_\i\leq C \|[u, R_iR_j](\r u)\|_p^{1-\f4p}\|\nabla\big([u, R_iR_j](\r u)\big)\|_{\f{4p}{p+4}}^{\f4p}\\
 \di \leq C \Big[\|u\|_{\rm BMO}\|\r u\|_p\Big]^{1-\f4p}\Big[\|\nabla u\|_4\|\r u\|_{p}\Big]^{\f4p}\\[3mm]
 \di \leq C\|\nabla u\|_2^{1-\f4p}\|\nabla u\|_4^{\f4p}\|\r u\|_p \di \leq C\Phi_T^{1+\f\b4}\Big(\|\nabla u\|_2^{1-\f2p}+1\Big)\|\nabla u\|_2^{1-\f4p}\|\nabla u\|_4^{\f4p},
  \end{array}
 $$
 while
 $$
 \begin{array}{ll}
 \di \|\nabla u\|_4\leq C\big(\|{\rm div} u\|_4+\|\o\|_4\big)\leq C\big(\|\f{F+(P(\r)-P(\bar\r))}{2\mu+\l(\r)}\|_4+\|\r\|_\i^{\f14}\|\o\|_2^{\f12}\|\sqrt\r \dot u\|_2^{\f12}\big)\\
 \di \leq C\big(\|\f{F^2}{\sqrt{2\mu+\l(\r)}}\|_2^{\f12}+1+\|\r\|_\i^{\f14}\|\o\|_2^{\f12}\|\sqrt\r \dot u\|_2^{\f12}\big)\\
 \di \leq C\big(\|\f{F}{\sqrt{2\mu+\l(\r)}}\|_2^{\f{1-\v}{2}}\|\r\|_\i^{\f{1+\b\v}{4}}\|\f{F}{\sqrt{2\mu+\l(\r)}}\|_2^{\f\v2}\|\sqrt\r \dot u\|_2^{\f12}+1+\|\r\|_\i^{\f14}\|\o\|_2^{\f12}\|\sqrt\r \dot u\|_2^{\f12}\big)\\
 \di \leq C\big(\|\r\|_\i^{\f{1+\b\v}{4}}+1\big)\Big[\|\f{F}{\sqrt{2\mu+\l(\r)}}\|_2^{\f12}\|\sqrt\r \dot u\|_2^{\f12}+\|\o\|_2^{\f12}\|\sqrt\r \dot u\|_2^{\f12}+1\Big]\\
 \di \leq C\big(\|\r\|_\i^{\f{1+\b\v+\b}{4}}+1\big)(e+\|\nabla u\|_2)\Big(\f{\varphi^2(t)}{e+Z^2(t)}\Big)^{\f14},
 \end{array}
 $$
 where in the fourth inequality one has used the fact \eqref{ue3}.
 Then it holds that
 $$
  \begin{array}{ll}
 \di \|[u, R_iR_j](\r u)\|_\i\leq C\big(\|\r\|_\i^{1+\f\b4+\f{1+\b\v+\b}{p}}+1\big)(e+\|\nabla u\|^2_2)^{1-\f1p}\Big(\f{\varphi^2(t)}{e+Z^2(t)}\Big)^{\f1p}\\
 \di \leq C\Big(\f{\varphi^2(t)}{e+Z^2(t)}+1\Big)+C\big(\|\r\|_\i^{\big[1+\f\b4+\f{1+\b\v+\b}{p}\big]\f{p}{p-1}}+1\big)(e+\|\nabla u\|^2_2).
   \end{array}
 $$
 Thus it holds that for any $\v>0$, one can choose sufficiently large $p>2$ such that
\begin{equation}\label{eta1}
  \begin{array}{ll}
 \di \int_0^T \|[u, R_iR_j](\r u)\|_\i(t)dt \leq C\Phi_T^{1+\f\b4+\v}.
   \end{array}
\end{equation}
 By Lemma \ref{lemma5}, it holds that for suitably large but fixed $m>1$,
 $$
 \|\x_1+\x_2\|_{2m}\leq C m^\f12\Big[\|\r-\bar\r\|_{2m}+\|\sqrt\r u |x|^{\f\a2}\|_2^{\f{2}{m\a}}\Big]\leq C_m.
 $$
 Then
 $$
 \|\nabla (\x_1+\x_2)\|_2\leq C\|\r u\|_2\leq C\|\r\|_\i^{\f12}\|\sqrt\r u\|_2\leq C\|\r\|_\i^{\f12},
 $$
 and then
 $$
 {\rm log}^{\f12}(e+ \|\nabla (\x_1+\x_2)\|_{2m})\leq C ~{\rm log}^{\f12}(e+ \|\r u\|_{2m})\leq C_m ~{\rm log}^{\f12}(e+ \|\nabla u\|_{2})\leq C\Phi_T^{\f{1+\b\v}{2}}.
 $$
 Therefore, it holds that
\begin{equation}\label{x11}
 \|\x_1+\x_2\|_\i\leq C\big( \|\x_1+\x_2\|_{2m}+ \|\nabla (\x_1+\x_2)\|_2\big){\rm log}^{\f12}(e+ \|\nabla (\x_1+\x_2)\|_{2m})\leq C\Phi_T^{1+\f{\b\v}{2}}.
\end{equation}
 Finally, substituting \eqref{eta1} and \eqref{x11} into \eqref{theta-1}, it holds that
 $$
\Phi_T^\b\leq C \Phi_T^{1+\f\b4+\v}+C.
 $$
 Therefore, if $\b>\f43$ and choose $\v$ suitably small, then
\begin{equation}\label{up-d}
 \sup_{t\in[0,T]}\|\r\|_\i(t)\leq C_1,
\end{equation}
for some positive constant $C_1$. Again by \eqref{theta-1}, \eqref{eta1}, \eqref{x11} and \eqref{up-d}, it holds that
$$
\sup_{t\in[0,T]}\|\ln \r(t,\cdot)\|_\i\leq C,
$$
which implies that there exists some positive constant $c_1$ such that
  \begin{equation*}
    \r(t,x)\geq c_1>0,\qquad \forall (t,x)\in[0,T]\times \mathbb{R}^2.
    \end{equation*}
Thus the proof of Lemma \ref{lemma-den} is completed. $\hfill\Box$

\section{Proof of main results}

In this section, we give a sketch of proof of our main results.

 {\bf
Proof of Theorem \ref{theorem1}:}

Under the assumptions of the
theorem, the local existence of the classical solution can be
proved in a similar way as in \cite{Luo, S} and we omit it for
simplicity. In view of the lower and upper bound of the density
obtained in Section 3, the compressible Navier-Stokes equations
\eqref{CNS} are a hyperbolic-parabolic coupled system. One
can get the higher order a priori estimates.  Using these a priori
estimates, one can extend the local solution to the global one in a
standard way(see \cite{JWX2, JWX3} for more details). The proof of
Theorem \ref{theorem1} is complete.

{\bf Proof of Theorem \ref{theorem2}:}

To use Theorem \ref{theorem1}, we first construct the approximation of the initial data in \eqref{in-d} as follows.
 Since $\di\lim_{|x|\rightarrow+\i} \r_0(x)=\bar\r>0$,  there exists a large number $M>0$
 such that if $|x|\geq M$,  $\r_0(x)\geq \frac{\bar\rho}{2}.$
Then for any $0<\d<\f{\bar\r}{2}$, we define
\begin{equation}\label{app-in-d}
\r_0^\d(x)=\left\{
\begin{array}{ll}
\di\r_0(x)+\d,~~{\rm if}~|x|\leq M,\\
\di \r_0(x)+\d s(x),~~{\rm if}~M\leq|x|\leq M+1,\\
\di \r_0(x),~{\rm if}~ |x|\geq M+1,
\end{array}
\right.
\end{equation}
where $s(x)=s(|x|)$ is a smooth and decreasing function satisfying
$s(x)\equiv1$ if $|x|\leq M$ and $s(x)=0$ if $|x|\geq M+1$.
Similarly, one can construct the approximation of the initial
pressure denoted by $P_0^\d(x)$. Then it follows that $(\r_0^\d,
P_0^\d)(x)$ are regular functions satisfying $\r_0^\d(x)>\d,
P_0^\d(x)>P(\d)$ for any $x\in\mathbb{R}^2$ and
$(\r_0^\d,P_0^\d)(x)=(\r_0,P_0)(x)$ if $|x|\geq M+1$. Moreover, one
has
\begin{equation*}
(\r_0^\d-\bar\r, P_0^\d-P(\bar \r))\rightarrow (\r_0-\bar\r,
P(\r_0)-P(\bar\r)) ~~{\rm in}~W^{2,q}(\mathbb{R}^2)\times
W^{2,q}(\mathbb{R}^2),
\end{equation*}
and
\begin{equation*}
\Psi(\r_0^\d,\bar\r)(1+|x|^\a)\rightarrow
\Psi(\r_0,\bar\r)(1+|x|^\a) ~~{\rm in}~L^1(\mathbb{R}^2),
\end{equation*}
as $\d\rightarrow 0$. To construct the approximation of the initial
velocity, we define $u_0^\d$ as
\begin{equation}\label{new-cc}
u_0^\d=\left\{
\begin{array}{ll}
 \di \tilde{u}_0^\d,\qquad |x|\leq M+1,\\
 \di u_0,\qquad |x|\geq M+1,
\end{array}
\right.
\end{equation}
where $ \tilde{u}_0^\d$ is the unique solution to the following
elliptic problem
\begin{equation}\label{n1}
\left\{
\begin{array}{ll}
 \di \mathcal{L}_{\r_0^\d}\tilde{u}_0^\d=\nabla P_0^\d+\sqrt{\r_0}g,\qquad {\rm in}~~\Omega_M:=\{x|\ |x|<M+1\},\\
 \di \tilde{u}_0^\d|_{|x|=M+1}=u_0.
\end{array}
\right.
\end{equation}
From \eqref{n1}, one has
\begin{equation}\label{new-cc1}
  \mathcal{L}_{\r_0} \tilde{u}_0^\d=-\nabla\big[(\l(\r_0^\d)-\l(\r_0)){\rm div} \tilde{u}_0^\d\big]+\nabla P_0^\d+\sqrt{\r_0} g, ~~{\rm in}~ \Omega_M.
\end{equation}
By the elliptic regularity, one has
\begin{equation}\label{ud1}
\begin{array}{ll}
\di \|\tilde{u}_0^\d\|_{H^2(\Omega_M)}\\
\di\leq C\Big[\|\l(\r_0^\d)-\l(\r_0)\|_\i\|\nabla({\rm div} \tilde{u}_0^\d)\|_2+\|\nabla(\l(\r_0^\d)-\l(\r_0))\|_\i\|{\rm div} u_0^\d\|_2+\|\nabla P_0^\d\|_2+\|\sqrt{\r_0} g\|_2\Big]\\
\di\leq C\Big[\d\|\nabla^2\tilde{u}_0^\d\|_{2}+\|\nabla P_0^\d\|_{2}+\|\sqrt{\r_0}\|_{L^\i(\mathbb{R}^2)} \|g\|_2\Big]\\
\leq C\Big[\d\|\nabla^2 \tilde{u}_0^\d\|_{2}+1\Big].
\end{array}
\end{equation}
where the generic positive constant $C$ is independent of $\d>0.$
Therefore, it follows from \eqref{ud1} that
\begin{equation}\label{h2}
\| \tilde{u}_0^\d\|_{H^2(\Omega_M)}\di\leq C
\end{equation}
where the positive constant $C$ is independent of $0<\d\ll1.$

From the compatibility conditions \eqref{cc}, \eqref{n1} and \eqref{new-cc1}, it
holds that
\begin{equation}\label{h5}
\left\{
\begin{array}{ll}
\di \mathcal{L}_{\r_0}(\tilde{u}_0^\d-u_0)
=-\nabla\big[(\l(\r_0^\d)-\l(\r_0)){\rm div} u_0^\d\big]+\nabla (P_0^\d-P_0):=\T^\d,\quad {\rm in}~~\Omega_M,\\
\di \di (\tilde{u}_0^\d-u_0)|_{|x|=M+1}=0.
\end{array}
\right.
\end{equation}
 It follows from
\eqref{app-in-d}, \eqref{h2} and \eqref{h5} that
\begin{equation}\label{h5+}
\tilde{u}_0^\d-u_0\in H_0^1(\Omega_M)\cap H^2(\Omega_M),
\end{equation}
and
\begin{equation}\label{h6}
\begin{array}{ll}
\|\tilde{u}_0^\d-u_0\|_{H^2(\Omega_M)}\leq C\|\T^\d\|_2\\
\di\leq\di C\Big[\|\l(\r_0^\d)-\l(\r_0)\|_{L^\i(\Omega_M)}\| \nabla^2\tilde{u}_0^\d\|_{2}+\|\nabla(\l(\r_0^\d)-\l(\r_0))\|_{L^\i(\Omega_M)}\| {\rm div} \tilde{u}_0^\d\|_{2}+\|\nabla(P_0^\d-P_0)\|_2\Big]\\
\di\leq
C\Big[\|\l(\r_0^\d)-\l(\r_0)\|_{L^\i(\Omega_M)}+\|\nabla(\l(\r_0^\d)-\l(\r_0))\|_{L^\i(\Omega_M)}+\|\nabla(P_0^\d-P_0)\|_2\Big]\\
\di\leq C\d~\rightarrow 0,
\end{array}
\end{equation}
as $\d\rightarrow 0$. It follows from \eqref{new-cc}, \eqref{h5+} and \eqref{h6} that
$u_0^\d\in H^2(\mathbb{R}^2)$ and
$$
u_0^\d\rightarrow u_0,~~{\rm in} ~ H^2(\mathbb{R}^2),
$$
and
$$
\sqrt{\r_0^\d} u_0^\d|x|^\a\rightarrow \sqrt\r_0u_0|x|^\a,~~{\rm in} ~ L^2(\mathbb{R}^2),
$$
as $\d\rightarrow0.$
 By Theorem \ref{theorem1},  there exists a
unique classical solution $(\r^\d,u^\d)$ to the compressible
Navier-Stokes equations \eqref{CNS} with the initial data $(\r_0^\d,
P_0^\d, u_0^\d)$ such that $c_{\d}\leq\r^\d\leq C$ for some positive
constants $c_\d$ depending on $\d$ and $C>0$. It should be noted
that the estimates obtained in Section 3  are independent of
 the lower bound of the initial density $\r_0(x)$ except the lower bound of the density $\r(t,x)$ in Lemma \ref{lemma-den}.
Then we can pass the limit
$\d\rightarrow0$ to get the classical solution satisfying
\eqref{Jan-16-1}. It is referred to \cite{JWX3} for more details and
the proof of Theorem \ref{theorem2} is finished.

\end{document}